\newif\ifpictures
\picturestrue
\documentclass[12pt]{amsart}
\usepackage{amssymb, colordvi}
\usepackage{graphicx} 

\copyrightinfo{}{}
\newcounter{FNC}[page]
\def\fauxfootnote#1{{\addtocounter{FNC}{2}$^\fnsymbol{FNC}$%
     \let\thefootnote\relax\footnotetext{$^\fnsymbol{FNC}$#1}}}
\newcommand{\calB}{\mathcal{B}}
\newcommand{\calL}{\mathcal{L}}
\newcommand{\calP}{\mathcal{P}}
\newcommand{\calQ}{\mathcal{Q}}
\newcommand{\calS}{\mathcal{S}}
\newcommand{\calT}{\mathcal{T}}

\newcommand{\G}{\mathbb{G}}
\renewcommand{\P}{{\mathbb P}}
\newcommand{\C}{\mathbb{C}}
\newcommand{\R}{\mathbb{R}}
\newcommand{\Z}{{\mathbb Z}}
\newcommand{\QED}{\ifpictures\includegraphics[height=12pt]{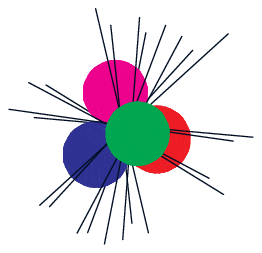}\else$\Box$\fi}
\newtheorem{thm}{Theorem}

\newtheorem{lemma}[thm]{Lemma}

\theoremstyle{definition}

\newtheorem{example}[thm]{Example}
\newtheorem{remark1}[thm]{Remark}

\newenvironment{ex}{\begin{example}\rm}{\end{example}}

\numberwithin{thm}{section}

\title{Line problems in nonlinear computational geometry}
\author{Frank Sottile}
\address{Department of Mathematics\\
         Texas A\&M University\\
         College Station\\
         Texas \ 77843\\
         USA}
\email{sottile@math.tamu.edu}
\urladdr{http://www.math.tamu.edu/~sottile/}

\author{Thorsten Theobald}
\address{Fachbereich Informatik und Mathematik\\
        J.W.\ Goethe-Universit\"at\\
        D--60054 Frankfurt am Main \\
        Germany}
\email{theobald@math.uni-frankfurt.de}
\urladdr{http://www.math.uni-frankfurt.de/~theobald/}

\thanks{Sottile was supported in part by NSF CAREER grant
  DMS-0538734 and Peter Gritzmann of the Technische
         Universit\"at M\"unchen.}
\thanks{Theobald was supported in part by
  DFG grant TH 1333/1-1.} 

\begin{document}

\begin{abstract}
 We first review some topics in the classical computational geometry of lines,
 in particular the $O(n^{3+\epsilon})$ bounds for the combinatorial complexity
 of the set of lines in $\R^3$ interacting with $n$ objects of fixed description
 complexity.
 The main part of this survey is recent work on a core algebraic
 problem---studying the lines tangent to $k$ spheres that also meet $4{-}k$
 fixed lines. 
\end{abstract}

\maketitle

%%%%%%%%%%%%%%%%%%%%%%%%%%%%%%%%%%%%%%%%%%%%%%%%%%%%%%%%%%%%%%%%%%%%%%%%%%%%%%
\section{Introduction}

While classical computational geometry is often concerned with (piecewise) linear
objects such as polyhedra, current applications and interests include the study
of  algorithms and computations with curved (non-linear)
objects~\cite{bhmmvt-2006,mishra-handbook,theobald-habil}.
In contrast to classical computational geometry, whose algorithms often have
low degree polynomial complexity,
this \Blue{{\it nonlinear computational geometry}} is often intractable as
many natural problems have intrinsic single- or double-exponential
complexity~\cite{bpr-03}. 

Algorithmic questions involving lines in $\R^3$ and $\R^d$ (which 
are fundamental problems in computational 
geometry \cite{cegss-96, pellegrini-handbook,stolfi-b91}) 
are at the boundary between classical and nonlinear
computational geometry.
Although a line in $\R^3$ is a \Blue{\emph{polyhedral}} set in $\R^3$
and collections of lines have a rich combinatorial structure,
the interaction of lines with other objects is inherently 
nonlinear, since the space of lines constitutes a curved 
submanifold (Grassmannian) in natural global Pl\"ucker coordinates.
However, many line problems have low algebraic degree and dimension
and are therefore computationally tractable, with solutions involving
both combinatorial and nonlinear techniques.

The investigation of line problems in computational geometry
arose from applications such as hidden surface removal, motion planning, and
configurations of mechanisms. 
Early work typically studied interactions of lines with lines and other 
polyhedral objects (see \cite{cegss-96,pellegrini-shor-92}).
In the last decade, several series of papers have studied problems 
of line transversals to more general objects (such as spheres) from various
viewpoints. From the combinatorial point of view, tight bounds on the combinatorial
complexity were obtained, while
nearly complete solutions to the underlying (real) algebraic geometric problems
were given.

We survey these developments in the computational geometry of lines.
Our particular goal is to relate the different viewpoints and to
explain some techniques (e.g., from computational algebraic
geometry) which are not common in discrete and computational geometry, 
but which have been 
fruitful in the study of line problems. 
These techniques provide useful tools for other 
problems in nonlinear computational geometry. We also exhibit 
geometric configurations which prove that certain bounds on the number of
real solutions are tight. Some of these constructions 
are new, in particular a configuration with four disjoint
spheres in $\R^3$ having 12 distinct common tangent lines and 
6 geometric permutations (see Figure~\ref{F:disjoint}).

The paper is structured as follows.
Section~\ref{S:classical} 
presents the classical computational geometry of lines in $\R^3$,
including the fundamental bound of $O(n^{3+\epsilon})$ for the combinatorial
complexity of the set of lines interacting with $n$ objects in $\R^3$.
The point of such bounds is that naive arguments give a bound of $O(n^4)$.
Section~\ref{se:algebraic} contains the heart of this survey, where we study the
algebraic core problem of lines tangent to $k$
spheres and transversal to $4{-}k$ lines in $\R^3$.
In Section~\ref{S:openproblems} we discuss some open problems.

%%%%%%%%%%%%%%%%%%%%%%%%%%%%%%%%%%%%%%%%%%%%%%%%%%%%%%%%%%%%%%%%%%%%%%%
%
%
%
\section{Classical computational geometry of lines}\label{S:classical}

%%%%%%%%%%%%%%%%%%%%%%%%%%%%%%%%%%%%%%%%%%%%%%%%%%%%%%%%%%%%%%%%%%%%%%%
%
%
%
\subsection{Pl\"ucker coordinates for lines}\label{S:Pluecker}
A source of non-linearity in the computational geometry of lines 
is that while lines are objects of linear algebra, the set of all lines is
naturally a curved submanifold of projective space.
If we represent a line in (projective) $d$-space $\P^d$ as the affine span of
two points $x^T=(x_0,x_1,\dotsc,x_d)$
and $y^T=(y_0,y_1,\dotsc,y_d)$, then its 
\Blue{{\it Pl\"ucker coordinates}}~\cite{pottmann-wallner-b2001} are
\[
   \Blue{p_{ij}}\ :=\ x_iy_j-x_jy_i\ =\ 
    \left|\begin{matrix}x_i&x_j\\y_i&y_j\end{matrix}\right|
   \qquad\mbox{for}\quad 0\leq i<j\leq d\,,
\]
which define a point in $D$-dimensional projective space $\P^D$, where 
$D := \binom{d+1}{2}-1$.

The set \Blue{$\G_{1,d}$} of all lines in $\P^d$ is called the
\Blue{\emph{Grassmannian of lines}} in $\P^d$.
The necessary and sufficient conditions for a point 
$(\,p_{ij}\mid0\leq i<j\leq d\,)\in\P^D$ to represent a line and hence lie in
$\G_{1,d}$ are furnished by the 
quadratic \Blue{\emph{Pl\"ucker equations}},
 \[   
   p_{ij} p_{kl} - p_{ik} p_{jl} + p_{il} p_{jk}\ =\ 0\qquad\mbox{for}\quad
    0\leq i<j<k<l\leq d\,.
 \] 
When $d=3$, the Grassmannian $\G_{1,3}$ is the hypersurface in $\P^5$ 
cut out by the single equation
 \begin{equation}\label{Eq:Pluecker}
   p_{01}p_{23}\ -\ p_{02}p_{13}\ +\ p_{03}p_{12}\ =\ 0\,.
 \end{equation}

Geometric conditions
such as incidence or tangency are naturally expressed in terms of 
Pl\"ucker coordinates.
For example, two lines in $\P^3$ meet if and only if their Pl\"ucker coordinates
$\Blue{(p_{ij})}$ and $\Magenta{(p'_{ij})}$ satisfy the bilinear equation
 \begin{equation}\label{Eq:linear_Pluecker}
  \Blue{p_{01}}\Magenta{p'_{23}}\ -\ 
  \Blue{p_{02}}\Magenta{p'_{13}}\ +\ 
  \Blue{p_{03}}\Magenta{p'_{12}}\ +\ 
  \Blue{p_{12}}\Magenta{p'_{03}}\ -\ 
  \Blue{p_{13}}\Magenta{p'_{02}}\ +\ 
  \Blue{p_{23}}\Magenta{p'_{01}}\ =\ 0\,.
 \end{equation}
Indeed, two lines spanned by points \Blue{$x,y$} and 
\Magenta{$x',y'$} in $\P^3$ meet if and only if the four points are affinely
dependent, which is expressed by $\det(\Blue{x,y},\Magenta{x',y'})=0$.
Laplace expansion of this determinant along the first two columns gives the
linear equation~\eqref{Eq:linear_Pluecker}.

A sphere in $\R^3$ with radius $r$ and center $c$ has equation 
$x^TQx=0$, where 
\[
   x^T\ =\ (1,x_1,x_2,x_3) \qquad\mbox{and}\qquad
   Q\ =\ \left(
    \begin{matrix}
      c_1^2 + c_2^2 + c_3^2 - r^2 & -c_1 & -c_2 & -c_3 \\
      -c_1 & 1 & 0 & 0 \\
      -c_2 & 0 & 1 & 0 \\
      -c_3 & 0 & 0 & 1
    \end{matrix}
   \right)  \,.
\]
A line $\ell$ spanned by two points $x$ and $y$ in $\P^3$ corresponds to a
2-dimensional linear subspace $H$ in $\R^4$ spanned by the vectors $x$ and $y$.
We may also consider $H$ to be a matrix with columns $x$ and $y$.

The line $\ell$ is tangent to the sphere if and only if the restriction of
the quadratic form $Q$ to $H$ is singular.
This restriction is represented by the $2\times 2$ symmetric matrix $H^T Q H$.
Taking its determinant, we obtain
 \[ 
   \bigl(\wedge^2H\bigr)^T\, ( \wedge^2Q )\, \wedge^2\!H\ =\ 0\,,
 \] 
where $\wedge^2 H$ is the vector of Pl\"ucker coordinates for the line $\ell$
and $\wedge^2 Q$ is the $6\times 6$ matrix whose entries are the $2\times 2$ 
minors of $Q$.
This is a quadratic equation in the Pl\"ucker coordinates of $\ell$ which cuts out
the lines tangent to the sphere defined by $Q$.\smallskip

A key tool for us is the following estimate.\medskip

\noindent{\bf B\'ezout's Theorem. }
{\it
 The number of isolated solutions to $n$ polynomial equations in $n$-space is 
 bounded from above by the product of the degrees of the polynomials.
}

%%%%%%%%%%%%%%%%%%%%%%%%%%%%%%%%%%%%%%%%%%%%%%%%%%%%%%%%%%%%%%%%%%%%%%%
%
%
%
\subsection{Combinatorial complexity}

We consider the combinatorial complexity of the set $\calT(\calS)$ of lines 
interacting (in a specified way) with a given set $\calS$ of lines or other
objects in $\R^3$.  
If the sets in $\calS$ are \Blue{{\it semi-algebraic}} in that they are defined by
equations and inequalities, then $\calT(\calS)$ will be a semi-algebraic set in
$\G_{1,3}$.
Its boundary $\partial\calT(\calS)$ is also semi-algebraic
and consists of objects which are tangent to at least one set in $\calS$.

For example, fix a line $\ell \subset \R^3$ 
and let $\calT(\ell)$ be the
set of lines which pass above $\ell$.
The boundary $\partial\calT(\ell)$ consists of lines that meet $\ell$,
and we have already seen that this boundary is defined by a linear equation in
the Pl\"ucker coordinates~\eqref{Eq:linear_Pluecker}.

For another example, 
let $\calT(C)$ be the set of lines which intersect a fixed convex body~$C$.
Its boundary $\partial\calT(C)$ consists of lines which are tangent to $C$.
When $C$ is a ball in $\R^3$ whose boundary sphere is given by an 
equation of the form $x^T Q x=0$, then  
$\calT(C)$ and $\partial\calT(C)$ are defined by the conditions
 \begin{equation}
  \label{eq:quadrictangentcond}
  p^T \, \bigl(\wedge^2 Q\bigr) \, p \ \leq \ 0 \qquad\text{and}\qquad
%
%   The \leq is really necessary here.  Frank Checked this twice.
%
  p^T \, \bigl(\wedge^2 Q\bigr) \, p \ = \ 0 \, ,
 \end{equation}
respectively.

A \Blue{\emph{face}} of $\calT(\calS)$ 
is a connected component of the set of lines in $\partial\calT(\calS)$ 
which are tangent to a fixed subset of $\mathcal{S}$. 
The \Blue{\emph{combinatorial complexity}} of the set $\calT(\calS)$ 
is the total number of its faces.
In this combinatorial analysis, we assume that the objects are in general
position (e.g.~every four balls have 12 common tangents, see Section~\ref{se:algebraic}).

Suppose that we have 
$n$ lines $\calL=\{\ell_1, \ldots, \ell_n\}\subset\R^3$.
The \Blue{\emph{upper envelope of $\calL$}} is the set of lines which pass above
all elements of $\calL$.
The following result was proved in~\cite{cegss-96}.

%%%%%%%%%%%%%%%%%%%%%%%%%%%%%%%%
\begin{thm}\label{T:theta_n^3}
 The maximum combinatorial complexity of the entire upper envelope of $n$
 lines in space is $\Theta(n^3)$.
\end{thm}
%%%%%%%%%%%%%%%%%%%%%%%%%%%%%%%%%

The significance of this cubic upper bound is that every four lines could have two
common transversals.
This gives possibly $O(n^4)$ zero-dimensional faces of the upper envelope.
Only those transversals which pass above the other lines in $\calL$ are
zero-dimensional faces, and Theorem~\ref{T:theta_n^3} says that
there are relatively few of these.

\begin{figure}[ht]
\ifpictures
\[
  \begin{picture}(270,180)(0,5)
   \put(0,0){\includegraphics[height=180pt]{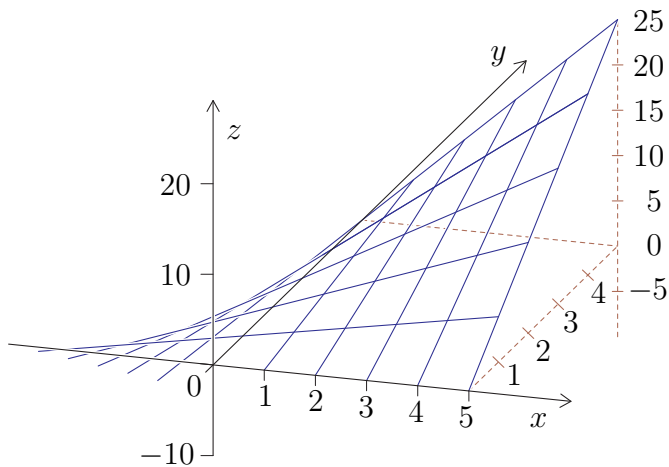}}
   \put( 60,33){$0$} \put( 87,29){$1$} \put(106,27){$2$}
   \put(125,25){$3$} \put(144,23){$4$} \put(163,21){$5$}
   \put(190,21){$x$}

   \put(75,130){$z$}
   \put(227, 69.5){$-5$} \put(234,86){$0$} \put(234,103){$5$} \put(229,120){$10$}
   \put(229,137){$15$} \put(229,154){$20$}\put(229,171){$25$}
   \put(42,  7){$-10$} \put(50,75){$10$} \put(50,109){$20$}

   \put(175,160){$y$} 
   \put(181,37){$1$}  \put(192,46){$2$}
   \put(202.5,57){$3$}\put(213,68){$4$}
  \end{picture}
\]
\fi

\caption{Families $\mathcal{G}$ and $\mathcal{H}$ of lines on the 
  hyperbolic paraboloid $z=xy$.\label{fi:hypparaboloid}}

\end{figure}

Here, $\Theta$ expresses that there is also a cubic lower bound, 
that is, there exists a construction whose upper envelope has exactly cubic
combinatorial complexity. 
Consider three families $\mathcal{G}$, $\mathcal{H}$, and $\calL$ of lines,
each of cardinality $N:=\lfloor n/3 \rfloor$. 
The lines of the
two families $\mathcal{G}$ and $\mathcal{H}$ are
\[
\begin{array}{rcl@{\quad}l}
  g_i & = & (0,i,0)^T + \R (1,0,i)^T \, , & 1 \le i \le N \, , \\
  h_j & = & (j,0,0)^T + \R (0,1,j)^T \, , & 1 \le j \le N \, , 
\end{array}
\]
which form a grid of lines on the hyperbolic
paraboloid $z = xy$.
Each line of $\mathcal{G}$ is parallel to the $xz$-plane
and each line of $\mathcal{H}$ is parallel to the
$yz$-plane (see Figure~\ref{fi:hypparaboloid}).
%%%%%%%%%%%%%%%%%%%%%%%%%%%%%%%%%%%%%%%%%%%%%%%%%%%%%%%%%%%%%%%%%%%%%%%%%%%%%%%%%
%%%%%%%%%%%%%%%%%%%%%%%%%%%%%%%%%%%%%%%%%%%%%%%%%%%%%%%%%%%%%%%%%%%%%%%%%%%%%
The third family of lines $\calL$ is given by
\[
  \ell_k \ = \ (0, 0, -n^5)^T + \R \Big(1,\frac{k}{n^2}-1,n \Big)^T , 
     \quad 1 \le k \le N \, .
\]
All lines in $\calL$ pass through the point $(0, 0, -n^5)^T$ below the
hyperboloic paraboloid and have a steep $z$ slope. 
For each triple of lines 
$(g_i, h_j, \ell_k) \in \mathcal{G} \times \mathcal{H}\times \calL$, 
there exists a line connecting some point of $\ell_k$ to the intersection point
$g_i \cap h_j$ which lies above all the other lines in
$\mathcal{G},\mathcal{H},\mathcal{L}$ (see \cite{cegss-96} for
details). 
This implies the lower bound $\Omega(n^3)$, and a
perturbation of this construction brings it into general position.

%%%%%%%%%%%%%%%%%%%%%%%%%%%%%%%%%%%%%%%%%%%%%%%%%%%%%%%%%%%%%%%%%%%%%%%%%%%%%%%%
%
%
%
\subsection{Line transversals to balls and semialgebraic sets}
Let $\calB = \{B_1, \ldots, B_n\}$ be a set of balls and $\calT(\calB)$ be the
set of lines which meet every ball in $\calB$.
Suppose that the balls are in general position (which is discussed in
Section~\ref{se:algebraic}).
Then, for each $j=0,1,2,3$, a $j$-dimensional face of $\calT(\calB)$
is a connected component of the set of lines 
in $\calT(\calB)$ which are tangent to a fixed set of $4{-}j$
balls in $\calB$. Agarwal, Aronov, and Sharir \cite{aas-99} have shown a
near-cubic bound on the complexity of $\calT(\calB)$. 
This result is striking---again there is a natural quartic upper bound.

\begin{thm} \label{th:combcomplexityballs}
 Let $\calB$ be a set of $n$ balls in $\R^3$. Then the complexity of
 $\calT(\calB)$ is $O(n^{3+\varepsilon})$ for any $\varepsilon > 0$.
\end{thm}

Note that the constant in the $O$-notation depends on $\varepsilon$.
For balls of general radius, the upper bound is essentially tight,
as there is a construction with complexity $\Omega(n^3)$.
For unit balls the tightness of the lower bound is open.
The best known construction has only quadratic complexity~\cite{aas-99}.

Theorem~\ref{th:combcomplexityballs} is now understood as a special case of a
more general result.
The \Blue{\emph{lower envelope}} of a collection of functions is their
pointwise minimum, and the  \Blue{\emph{upper envelope}} is their pointwise
maximum. 
A  \Blue{\emph{sandwich region}} is a region which lies above an upper
envelope for one collection of functions and below a lower envelope for a
different collection of functions.

The set of line transversals to $\calB$ is a sandwich region given by two sets
of trivariate functions (see~\cite{aas-99,cegss-96}). 
If we exclude lines parallel to the $yz$-plane (which is no loss of generality
as the balls are in general position), a line
$\ell$ in $\R^3$ can be uniquely represented by the quadruple
$(\sigma_1, \sigma_2, \sigma_3, \sigma_4) \in \R^4$
parametrizing its projections to the $xy$-
and $xz$-planes: $y = \sigma_1 x + \sigma_2$
and $z = \sigma_3 x + \sigma_4$.

Let $B$ be a ball in $\R^3$. For fixed
$\sigma_1,\sigma_2,\sigma_3$, the set of lines
$(\sigma_1,\sigma_2,\sigma_3,\sigma_4)$ that intersect $B$
is obtained by translating a line in the $z$-direction
between two extreme values
$(\sigma_1,\sigma_2,\sigma_3,\phi_B^{-}(\sigma_1,\sigma_2,\sigma_3))$ and
$(\sigma_1,\sigma_2,\sigma_3,\phi_B^{+}(\sigma_1,\sigma_2,\sigma_3))$,
which represent lines tangent to $B$ from below and from above,
respectively. Hence, the set of line transversals to $\calB$
is
\[
  \left\{(\sigma_1,\sigma_2,\sigma_3,\sigma_4) \: : \:
  \max_{B \in \calB} \phi_B^{-}(\sigma_1,\sigma_2,\sigma_3)
  \le \sigma_4 \le
  \min_{B \in \calB} \phi_B^{+}(\sigma_1,\sigma_2,\sigma_3) \right\} ,
\]
which is a sandwich region of two sets of trivariate functions. 
As the balls in $\calB$ are in general position, the vertices 
of the boundary of this region are  lines which are tangent
to four of the balls in $\calB$.

Recently, Koltun and Sharir have shown the following result for
trivariate  functions \cite{koltun-sharir-2003}.
Together with the preceeding discussion, this 
implies Theorem~\ref{th:combcomplexityballs}.

%%%%%%%%%%%%%%%%%%%%%%%%%%%%%%%%%%%%%%%%%%%%%%%%%%%%%%%%%%%%%%%%%%%%%
\begin{thm} \label{th:arrangements1}
 For any $\varepsilon > 0$, the combinatorial complexity
 of the sandwich region between two families of $n$ trivariate functions of
 constant description complexity is $O(n^{3+\varepsilon})$. 
\end{thm}
%%%%%%%%%%%%%%%%%%%%%%%%%%%%%%%%%%%%%%%%%%%%%%%%%%%%%%%%%%%%%%%%%%%%%

This is a central achievement in a sequence of papers
(including \cite{halperin-sharir-94, koltun-sharir-2003,sharir-94}). 
The basic idea is to reduce the geometric problem to
a combinatorial problem involving Davenport-Schinzel sequences.

Given a word $u=u_1 u_2\dotsc u_m$ over an alphabet with $n$ symbols, a
subsequence of $u$ is a word $u_{i_1} \ldots u_{i_k}$ for some indices
$1 \le i_1 < \cdots < i_k \le m$.
An \Blue{$(n,s)$-\emph{Davenport-Schinzel sequence}} is a word $u$ over an
alphabet with $n$ symbols such that
\begin{enumerate}
  \item two consecutive symbols of $u$ are distinct, and
  \item for any two symbols $a$, $b$ in the alphabet, the alternating
        sequence $abab\ldots\,$ of length $s+2$ is not a subsequence of $u$.
\end{enumerate}

\begin{figure}[ht]
\[
  \begin{picture}(300,150)
\ifpictures
   \put(5,0){ \includegraphics[height=150pt]{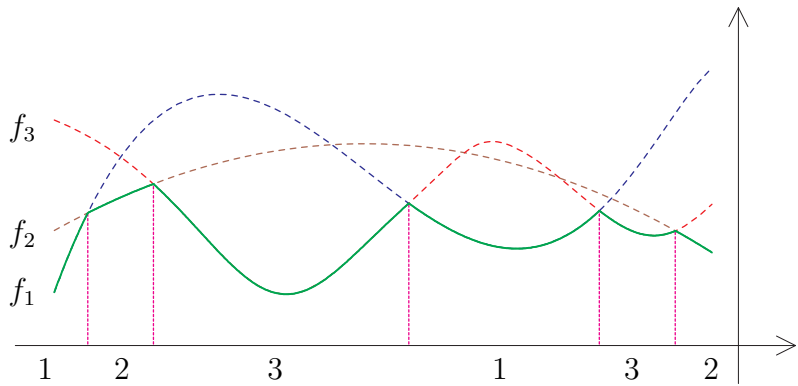}}
\fi
  \put(2,38){$f_1$}  \put(2,60){$f_2$}  \put(2,100){$f_3$}
  \put( 13,8){$1$}  \put( 42,8){$2$}\put(100,8){$3$}
  \put(185,8){$1$}  \put(235,8){$3$}\put(265,8){$2$}
  \end{picture}
\]
\caption{The lower envelope of three functions $f_1$, $f_2$, $f_3$ 
   with label sequence $(1,2,3,1,3,2)$.}
\label{fi:davenport}
\end{figure}

The lower envelope of $n$ continuous univariate functions $f_1,\dotsc,f_n$ 
is formed by a sequence of curved edges, where
each edge is a maximal connected subset of the envelope that belongs
to the graph of a single function $f_i$. 
This is illustrated in Figure~\ref{fi:davenport}. 
%%%%%%%%%%%%%%%%%%%%%%%%%%%%%%%%%%%%%%%%%%%%%%%%%%%%%%%%%%%%%%%%%%%%%%%%%%%%%
%%%%%%%%%%%%%%%%%%%%%%%%%%%%%%%%%%%%%%%%%%%%%%%%%%%%%%%%%%%%%%%%%%%%%%%%%%%

We label each edge by the index of the corresponding function.
Reading these labels left-to-right gives a sequence of labels.
If the graphs of the functions have
pairwise at most $s$ intersection points, then this sequence is an
$(n,s)$-Davenport-Schinzel sequence. 
Indeed, suppose that we have intervals $I$ and $J$ with distinct labels $i$ and
$j$ so that $f_i$ is the minimum of all functions over $I$ and $f_j$ is the
minimum of all functions over $J$. 
Then $f_i$ and $f_j$ must meet at some point that lies weakly between $I$ and
$J$.
Thus an alternating subsequence of length $s+2$ for two
symbols $i$ and $j$ implies the existence of $s+1$ intersection points between
the graphs of $f_i$ and $f_j$, which is a contradiction.
A (sharp) upper bound for the number of $(n,s)$-Davenport-Schinzel sequences is
only slightly super-linear in $n$ (rather than $n^2$). See,
e.g.,~\cite{sharir-agarwal-95}.\smallskip  

Recently, Agarwal, Aronov, Koltun, and Sharir~\cite{aaks-2005}
proved an upper bound on the 
combinatorial complexity of a set of lines defined with respect to a collection
of balls which does not fall in the framework of Theorem~\ref{th:arrangements1},
since it is not known whether this set can be realized as a 
sandwich region.
A line in $\R^3$ is \Blue{\emph{free}} (with respect to a set of balls $\calB$)
if it does not intersect the interior of any ball in $\calB$.

%%%%%%%%%%%%%%%%%%%%%%%%%%%%
\begin{thm}
For any $\varepsilon > 0$
the combinatorial complexity of the space of lines free with respect to a set of
$n$ unit balls in $\R^3$ is $O(n^{3+\varepsilon})$.
\end{thm}
%%%%%%%%%%%%%%%%%%%%%%%%%%%%

For polyhedra, the following result is known \cite{pellegrini-shor-92}.
Here, the \Blue{\emph{total complexity}} of a set $\calP$ of polyhedra 
is the total number of their faces.
Recall that $\calT(\calP)$ is the set of lines which meet every polyhedron in $\calP$.

%%%%%%%%%%%%%%%%%%%%%%%%%%%%%%%%%%%%%%%%%%%%%%%%%%%%%%%%
\begin{thm} Given a set $\mathcal{P}$ of polyhedra with total complexity $n$,
the combinatorial complexity of $\calT(\calP)$ is bounded by
$O(n^3 2^{c \sqrt{\log n}})$, where $c$ is a constant. Moreover, 
the set of extremal stabbing lines (those representing the vertices
of $\mathcal{T}(\calP)$) can be found in time $O(n^3 2^{c \sqrt{\log n}})$. 
\end{thm}
%%%%%%%%%%%%%%%%%%%%%%%%%%%%%%%%%%%%%%%%%%%%%%%%%%%%%%%%

In a similar spirit, in \cite{bdd-2004} it was shown
that a set $\mathcal{P}$ of $k$ polyhedra with total complexity~$n$
has $O(k^2 n^2)$ many common tangent lines. For $k = o(n)$ this gives
another example where the natural upper bound of $O(n^4)$ is not sharp. 

%%%%%%%%%%%%%%%%%%%%%%%%%%%%%%%%%%%%%%%%%%%%%%%%%%%%%%%%%%%%%%%%%%%%%%%%%%%%%%%%%
%
%
%
\subsection{Transversals, convexity, and geometric permutations}

A line transversal to a set $\mathcal{S}$ of pairwise disjoint convex bodies
in $\R^d$ induces two linear orders on $\mathcal{S}$. 
These two orders are reverse to each other, and we consider them as a single
\Blue{\emph{geometric permutation}}. 
Studying the geometry of line transversals and geometric
permutations
has a long history in discrete and computational geometry; see the
surveys~\cite{gpw-93,wenger-2004}. 
An obstacle
is the lack of a suitable convexity structure on the space of lines.
(Goodman and Pollack~\cite{GP95} discuss limitations of any notion of convexity on
Grassmannians.)
We discuss a recent result sowing that the set of directions to lines
realizing a given geometric permutation is convex.

Let $\mathcal{S}$ be an (ordered) sequence of convex bodies. 
A \Blue{\emph{directed line transversal}} 
to $\mathcal{S}$ is an oriented line intersecting all bodies
in the order on $\mathcal{S}$.
Let $\mathcal{K}(\mathcal{S})$ denote the set of  
directions of the directed line transversals to $\mathcal{S}$.

%%%%%%%%%%%%%%%%%%%%%%%%%%%%
\begin{thm} \label{th:inflatable}
 Let $\calB$ be a sequence of disjoint balls in $\R^d$. If 
 $\mathcal{K}(\calB)$ is nonempty, then it is strictly convex.
\end{thm}
%%%%%%%%%%%%%%%%%%%%%%%%%%%%

This result is due to Borcea, Goaoc, and Petitjean~\cite{bgp-2007}.
The main case involves three balls in $\R^3$.
They use a result from~\cite{cghp-2006} about the boundary of the set of directions to the family of
lines realizing a given geometric permutation.
From there, they study the convexity of this boundary curve to obtain their result.

To show convexity when $|\calB| > 3$, it suffices to 
show $\mathcal{K}(B) = 
\bigcap_{\calB' \subset \calB, | \calB'| = 3} \mathcal{K}(\mathcal{B'})$. 
The inclusion $\subset$ is obvious;
for the converse direction let
$v \in \bigcap_{\calB' \subset \calB, | \calB'| = 3} \mathcal{K}(\calB')$.
By choice of $v$, the orthogonal projections of any three balls 
along $v$ intersect. Hence,
the classical Helly theorem implies that the projections of all balls have a common 
point of intersection, and so the balls of $\calB$
have a common line transversal with direction $v$. Since this
transversal is consistent with the ordering induced by $\calB$ on
every triple of balls, it is a directed line transversal to $\calB$.

For convexity in $\R^d$ for $d>3$, two transversals with the
same geometric permutation lie in a 3-dimensional subspace.
Thus convexity in $\R^3$ implies convexity in $\R^d$.

This result caps a long series of papers.
First, Hadwiger~\cite{hadwiger-1957} proved this for collections of balls in which
the centers of any two balls are further apart than twice the sum of their radii.
Then Holmsen, Katchalski, and Lewis~\cite{holmsen-katchalski-lewis-2003} proved this
for disjoint unit balls, and it was extended Cheong, Goaoc, Holmsen, and Petitjean
to pairwise inflatable balls~\cite{cghp-2006}.
A set $\calB$ of balls with centers $c_i$ and radii $r_i$
is called  \Blue{\emph{pairwise inflatable}} if for every two balls 
$B_i, B_j \in \calB$ we have $||c_i-c_j||^2 > 2(r_i^2 + r_j^2)$, 
where $||\cdot||$ denotes the Euclidean norm.
In particular, any set of disjoint unit balls is pairwise inflatable.

%%%%%%%%%%%%%%%%%%%%%%%%%%%%%%%%%%%%%%%%%%%%%%%%%%%%%%%%%%%%%%%%%%%%%%%%%%%%%%
%
%
%
\section{The Algebraic Core\label{se:algebraic}}

Many problems from Section~\ref{S:classical}
share the following algebraic core.
%%%%%%%%%%%%%%%%%%%%%%%%%%%%%%%%%%%%%%%%%%%%%%%%%%%%%%%%%%%%%%%%%%%%%%%%%%%%
 \begin{equation}\label{Eq:Alg_core}
  \raisebox{6pt}{\begin{minipage}[t]{330pt}
   Given $k$ spheres and $4{-}k$ lines in $\R^3$, determine
   which common tangents to the spheres also meet each of $4{-}k$ given lines.
  \end{minipage}}
 \end{equation}
%%%%%%%%%%%%%%%%%%%%%%%%%%%%%%%%%%%%%%%%%%%%%%%%%%%%%%%%%%%%%%%%%%%%%%%%%%%
This basic question about the geometry of lines is also motivated by problems of
visibility in $\R^3$.
When viewing a scene in a particular direction along a moving viewpoint, the
objects that can be viewed may change when the line of sight becomes tangent to
one or more objects.
If the objects are in general position, then the most degenerate such lines of
sight are the lines tangent to four objects.
When the objects are polyhedra and spheres, these most 
degenerate lines of sight are the tangents to $k$ spheres which also meet $4{-}k$ edges.
Our algebraic core~\eqref{Eq:Alg_core} is the relaxation where we replace the
edges by their supporting lines.

We will discuss the number of lines (both real and complex) that can be
solutions to~\eqref{Eq:Alg_core} when the spheres and lines are in
general position.
By general position, we mean in the algebraic sense:
Every complex line solving~\eqref{Eq:Alg_core} is a simple solution (it occurs without
multiplicity), lines solving three 
of the conditions form a smooth 1-dimensional family, etc.
We will also classify the degenerate positions of $k$ spheres and $4{-}k$ lines which
admit a 1-dimensional family of solutions to~\eqref{Eq:Alg_core}.
Identifying these degeneracies is important when devising robust algorithms for
problems which involve this algebraic core.
Both of these investigations are also interesting from the point of view of
enumerative real algebraic geometry~\cite{ERAG} and of computational algebraic
geometry.

%%%%%%%%%%%%%%%%%%%%%%%%%%%%%%%%%%%%%%%%%%%%%%%%%%%%%%%%%%%%%%%%%%%
\subsection{The problem of four lines}\label{SS:4lines}

Let us begin with the simplest case of our algebraic core problem,
the classical and surprisingly non-linear problem of common
transversals to four lines in space.

\begin{figure}[htb]
 \[
   \begin{picture}(235,115)(0,-12)
    \put(0,-12){\ifpictures\includegraphics[height=115pt]{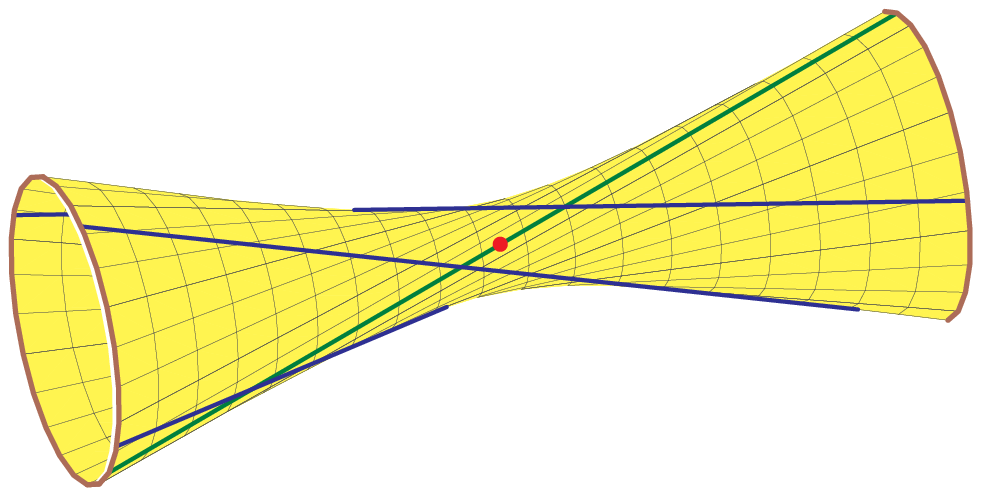}\fi}
    \put( 70, 5){$\ell_1$}
    \put(155,23){$\ell_2$}
    \put(100,60){$\ell_3$}
    \put(137,85){$p$}\put(136,81){\vector(-1,-2){16}}
    \put(190,104){$m_p$}
    \put(107,-12){$(i)$}
   \end{picture}
\qquad\quad
   \begin{picture}(170,105)(0,-12)
    \put(0,5){\ifpictures\includegraphics[height=90pt]{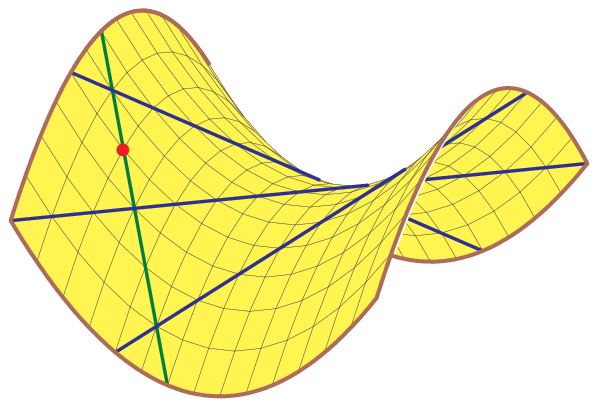}\fi}
    \put(120,25){$\ell_1$}
    \put(150,58){$\ell_2$}
    \put(135,80){$\ell_3$}
    \put( -6,79){$p$}\put( 1,77){\vector(2,-1){27}}
    \put( 13,99){$m_p$}
    \put(65,-12){$(ii)$}
   \end{picture}
\]
 \caption{Hyperboloids through 3 lines.\label{F:3lines}}
\end{figure}

Let $\ell_1$, $\ell_2$, $\ell_3$, and $\ell_4$ be lines in
general position in $\R^3$.
The three mutually skew lines $\ell_1$, $\ell_2$, and $\ell_3$ lie in one ruling
of a doubly-ruled hyperboloid as shown in Figure~\ref{F:3lines}.
%%%%%%%%%%%%%%%%%%%%%%%%%%%%%%%%%%%%%%%%%%%%%%%%%%%%%%%%%%%%%%%%%%%%%%%%%%%
%%%%%%%%%%%%%%%%%%%%%%%%%%%%%%%%%%%%%%%%%%%%%%%%%%%%%%%%%%%%%%%%%%%%%%%%%%%
This is either $(i)$ a hyperboloid of one sheet, or $(ii)$ a hyperbolic
paraboloid. 
The line transversals to $\ell_1$, $\ell_2$, and $\ell_3$ constitute the second
ruling. 
Through every point $p$ of the hyperboloid there is a unique line $m_p$ in the
second ruling which meets the lines $\ell_1$, $\ell_2$, and $\ell_3$.
To simplify this and other discussions, we will at times work in
projective space where parallel lines meet.
Thus for example, we include lines in the second ruling of $(i)$ which are
parallel to one of $\ell_1$, $\ell_2$, or $\ell_3$.

The hyperboloid is defined by a quadratic polynomial and so 
the fourth line $\ell_4$ will either meet the hyperboloid in two
points or  it will miss  the hyperboloid.
In the first case there will be two transversals to $\ell_1$, $\ell_2$, $\ell_3$,
and $\ell_4$, and in the second case there will be no transversals.
But in this second case $\ell_4$ meets the hyperboloid in two conjugate complex points, 
giving two complex transversals.

If $\ell_4$ is not in general position, it could be tangent to the hyperboloid,
in which case there is a single common transversal occurring with a multiplicity
of 2.
More interestingly, it could lie in the ruling of the first three
lines, and then there will be infinitely many common transversals.
\[
  \begin{picture}(245,125)
    \put(0,0){\ifpictures\includegraphics[height=120pt]{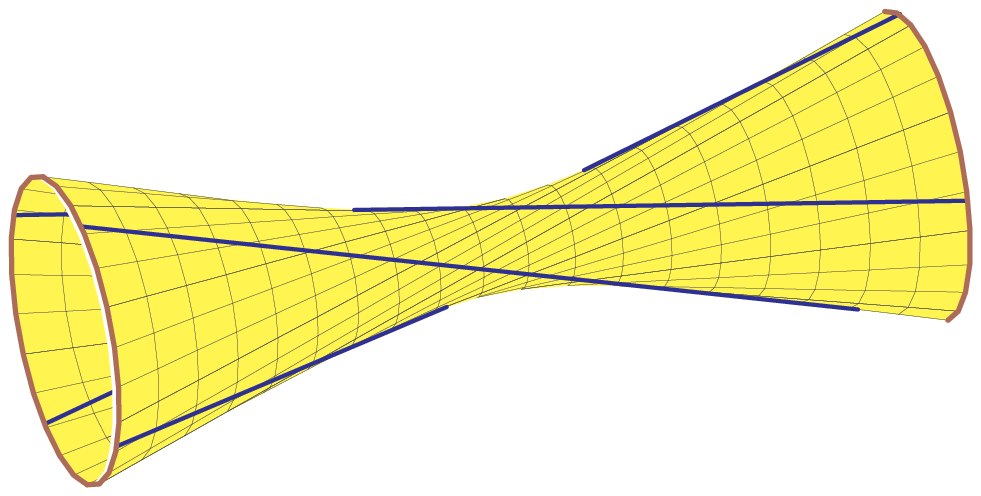}\fi}
    \put( 70, 15){$\ell_1$}
    \put(170, 35){$\ell_2$}
    \put( 95, 75){$\ell_3$}
    \put(180,110){$\ell_4$}
  \end{picture}
\]

Another possibility is for two of the first three lines to meet.
If $\ell_1$ and $\ell_2$ meet in a point $p$,
then they span a plane $H$, and the common transversals to $\ell_1$ and $\ell_2$
either lie in the plane $H$ or meet the point $p$.
In either case, further questions of transversals and tangents essentially
become problems in a plane.
Because of this reduction, we shall (almost always) assume that our lines are
pairwise skew. \smallskip

A different perspective on this problem of four lines is furnished by Pl\"ucker
coordinates for lines.
Recall from Section~\ref{S:Pluecker} that the set $\G_{1,3}$ of lines in $\P^3$
is a quadratic hypersurface in $\P^5$ defined by~\eqref{Eq:Pluecker}, and the
condition for a line $\ell$ to meet a fixed line is a linear
equation~\eqref{Eq:linear_Pluecker} in the Pl\"ucker coordinates of $\ell$.
Thus the Pl\"ucker coordinates of line transversals to 
$\ell_1$, $\ell_2$, $\ell_3$, and $\ell_4$ satisfy one quadratic and four linear
equations. 
By B\'ezout's Theorem, we see again
that there will be two line transversals to four given lines.
Observe that this purely algebraic derivation of the number of solutions cannot
distinguish between real and complex solutions.
This is both a strength (finding complex solutions is a relaxation of the
problem of finding real solutions) and a weakness of using algebra to study
configurations of objects in $\R^3$.

%%%%%%%%%%%%%%%%%%%%%%%%%%%%%%%%%%%%%%%%%%%%%%%%%%%%%%%%%%%%%%%%%%%%%%%%%%
\subsection{Lines tangent to four spheres}

As in Section~\ref{S:Pluecker}, a line is tangent to a sphere
with equation $x^TQx=0$ when its Pl\"ucker coordinates $p$ satisfy the quadratic
equation
 \begin{equation}\label{Eq:Sphere_eq}
   p^T\wedge^2\!Q\,\, p\ =\ 0\,.
 \end{equation}
It follows that the Pl\"ucker coordinates of common tangents to 4 spheres in
$\R^3$ are defined by 5 quadratic equations.
These are the Pl\"ucker equation~\eqref{Eq:Pluecker} and 4 quadratic
equations~\eqref{Eq:Sphere_eq} expressing tangency to each sphere.
B\'ezout's Theorem gives a bound of $2^5=32$ for the number of tangent lines.
Unlike the problem of lines in Section~\ref{SS:4lines}, this is not a sharp
upper bound, as not all solutions are isolated and geometrically meaningful.

Indeed, a sphere with center $(c_1,c_2,c_3)^T\in\R^3$ and radius $r$ is defined by
the equation
\[
  (x_1-c_1x_0)^2 +  (x_2-c_2x_0)^2 +  (x_3-c_3x_0)^2 \ =\ r^2 x_0^2\,.
\]
Every sphere in $\R^3$ contains the imaginary \Blue{{\it spherical conic at
infinity}},
 \begin{equation}\label{Eq:Spherical_Conic}
   \Blue{C}\ \colon\ x_0=0\quad\mbox{and}\quad 
   x_1^2+x_2^2+x_3^2=0\,.
 \end{equation}
In fact, spheres are exactly the quadrics that contain the spherical conic. 

Each (necessarily imaginary) line at infinity that is tangent to $C$ is tangent
to every sphere, and there is a one-dimensional family of such lines.
Thus the equations for lines tangent to four spheres not only define the 
tangents in $\R^3$, but also this one-dimensional excess component 
consisting of geometrically meaningless lines tangent to $C$.

We treat this excess component by defining it away.
Represent a line $\ell$ in $\R^3$ by a point $p\in\ell$ 
and a direction vector $v\in\R\P^2$.
No such line can lie at infinity, so we are avoiding the excess 
component of lines at infinity tangent to $C$.

%%%%%%%%%%%%%%%%%%%%%%%%%%%%%%%%%%%%%%%%%%%%%%%%%%%%%
\begin{lemma}\label{L:cubic_quartic}
 The set of direction vectors $v\in\R\P^2$ of lines tangent to four spheres with
 affinely independent centers consists of the common solutions to a cubic and a
 quartic equation on $\R\P^2$.
 Each direction vector gives one common tangent, and there are at most $12$ common
 tangents. 
\end{lemma}
%%%%%%%%%%%%%%%%%%%%%%%%%%%%%%%%%%%%%%%%%%%%%%%%%%%%%

\noindent{\it Proof.}
We follow the derivation of~\cite[Lemma~13]{theobald-enumerative}, which
generalizes the argument given in~\cite{mpt-2001} to spheres of arbitrary radii.
For vectors $x,y\in\R^3$, let $x\cdot y$ be their ordinary Euclidean dot
product and write $x^2$ for $x\cdot x$, which is $\|x\|^2$.

Fix $p$ to be the point of $\ell$ closest to the origin, so that 
 \begin{equation}\label{Eq:pdotv}
   p\cdot   v\ =\ 0\,.
 \end{equation}

The line $\ell$ is tangent to the sphere with radius $r$ centered at $c\in\R^3$
if its distance to $c$ is $r$,
 \[
   \left\|(c-p)\times  \frac{v}{\|v\|} \right\|\ =\ r\,.
 \]
Squaring and clearing the denominator $v^2$  gives
 \begin{equation}\label{Eq:tan_to_sphere}
  \left[(c-p)\times v\right]^2\ =\ r^2 v^2\,.
 \end{equation}
This formulation requires that $v^2\neq 0$.
A line with $v^2=0$ has a complex direction vector and it meets the  
spherical conic at infinity.

We assume that one sphere is centered at the origin with radius $r$,
while the other three have centers and radii $(c_i,r_i)$ for $i=1,2,3$.
The condition for the line to be tangent to the sphere centered at the origin is
 \begin{equation}\label{Eq:p2er2}
   p^2\ =\ r^2\,.
 \end{equation}
For the other spheres, we expand~\eqref{Eq:tan_to_sphere}, use vector product
identities, and the equations~\eqref{Eq:pdotv} and~\eqref{Eq:p2er2} to
obtain the vector equation
 \begin{equation}\label{Eq:p(v)}
   2v^2 \left(\begin{matrix}c_1^T\\c_2^T\\c_3^T\end{matrix}\right)
   \cdot p\ =\ 
   -\left(\begin{matrix}(c_1\cdot v)^2\\(c_2\cdot v)^2\\(c_3\cdot v)^2
          \end{matrix}\right)\ +\ 
    v^2 \left(\begin{matrix} c_1^2+r^2-r_1^2\\
           c_2^2+r^2-r_2^2\\c_3^2+r^2-r_3^2\end{matrix}\right)\ .
 \end{equation}

Now suppose that the spheres have affinely independent centers.
Then the matrix $(c_1,c_2,c_3)^T$ appearing in~\eqref{Eq:p(v)} is invertible.
Assuming $v^2\neq 0$, we may use~\eqref{Eq:p(v)} to write  $p$ as 
a quadratic function of $v$.
Substituting this expression into equations~\eqref{Eq:pdotv} and~\eqref{Eq:p2er2}, 
we obtain a cubic and a quartic equation for $v\in\R\P^2$.
\hfill\QED\medskip

B\'ezout's Theorem implies that there are at most $3\cdot 4=12$ isolated
solutions to these equations, and over $\C$ exactly 12 if they are generic.
The equations are however far from generic as they involve
only 13 parameters while the space of quartics has 14 parameters and the space
of cubics has 9 parameters.

%%%%%%%%%%%%%%%%%%%%%%%%%%%%%%%%%%%%%%%%%%%%%%%%%%%%%%%%%%%%%%%%%
\begin{ex}\label{Ex:cubic_quartic}
 Suppose that the spheres have equal radii, $r$, and have centers at the vertices
 of a regular tetrahedron with side length $2\sqrt{2}$,
 \[
   (2,2,0)^T\,,\qquad
   (2,0,2)^T\,,\qquad
   (0,2,2)^T\,,\quad\mbox{and}\quad
   (0,0,0)^T\,.
 \]
 In this symmetric case, the cubic factors into three linear factors.
 There are real common tangents only if 
 $\sqrt{2}\leq r\leq 3/2$, and exactly 12 when the inequality is strict.
 If $r=\sqrt{2}$, then the spheres are pairwise tangent and there are
 three common tangents, one for each pair of non-intersecting edges of the
 tetrahedron.
 Each tangent has algebraic multiplicity~4.
 If $r=3/2$, then there are six common tangents, each of multiplicity 2.
 The spheres meet pairwise in circles of radius $1/2$ lying in the plane
 equidistant from their centers.
 This plane also contains the centers of the other two spheres, as well as 
 one common tangent which is parallel to the edge between those centers.

 Figure~\ref{F:cu_qu} shows the cubic (which consists of three lines supporting
 the edges of an equilateral triangle) and the quartic, in an affine piece of 
 the set $\R\P^2$ of direction vectors.
 The vertices of the triangle are the standard coordinate directions $(1,0,0)^T$,
 $(0,1,0)^T$, and $(0,0,1)^T$.
%%%%%%%%%%%%%%%%%%%%%%%%%%%%%%%%%%%%%%%%%%%%%%%%%%%%%%%%%%%%%%%%%%%%%%%%%%%%%%%%%%
%
%     This was wrong
%
% The vertices of the triangle are
% the direction vectors  from the origin to the centers $c_1$, $c_2$, and $c_3$
% of the other three spheres. 
 The singular cases, $(i)$ when $r=\sqrt{2}$ and 
 $(ii)$ when $r=3/2$, are shown first, and then 
 $(iii)$ when $r=1.425$.
 The 12 points of intersection in this third case are visible in the expanded
 view in $(iii')$.
%%%%%%%%%%%%%%%%%%%%%%%%%%%%%%%%%%%%%%%%%%%%%%%%%%%%%%%%%%%%%%%%%%%%%%%%%%%%%%%
\begin{figure}[htb]
\[
  \begin{picture}(139,125)(2,-5)
   \put(0,3){\ifpictures\includegraphics[height=120pt]{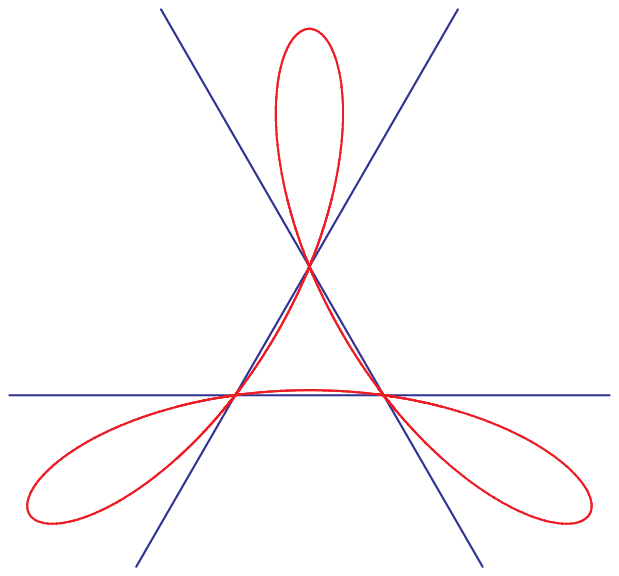}\fi}
%   \put( 52.5,118){$4$}
%   \put( 43  ,  9){$-2$}
%   \put( -3,56){$-3$} \put( 34,56){$-1$}
%   \put( 79,56){$1$}  \put(115,56){$3$}
   \put(56,0){$(i)$}
  \end{picture}
   \qquad\qquad
  \begin{picture}(135,125)(2,-5)
   \put(0,0){\ifpictures\includegraphics[height=120pt]{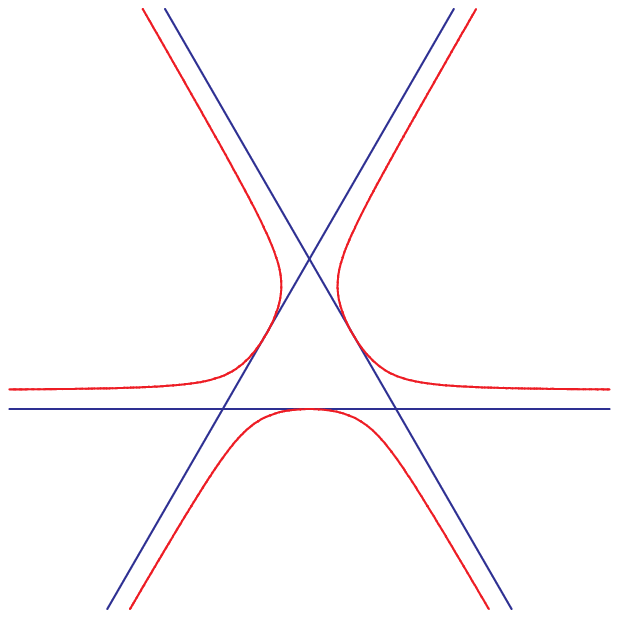}\fi}
%   \put( 52.5,109){$3$}
%   \put( 43,2.5){$-2$}
%   \put(  9,56){$-2$}
%   \put(103,56){$2$}  
   \put(54.5,0){$(ii)$}
  \end{picture}
\]
\[
  \begin{picture}(134,145)(1,-8)
   \put(0,1){\ifpictures\includegraphics[height=120pt]{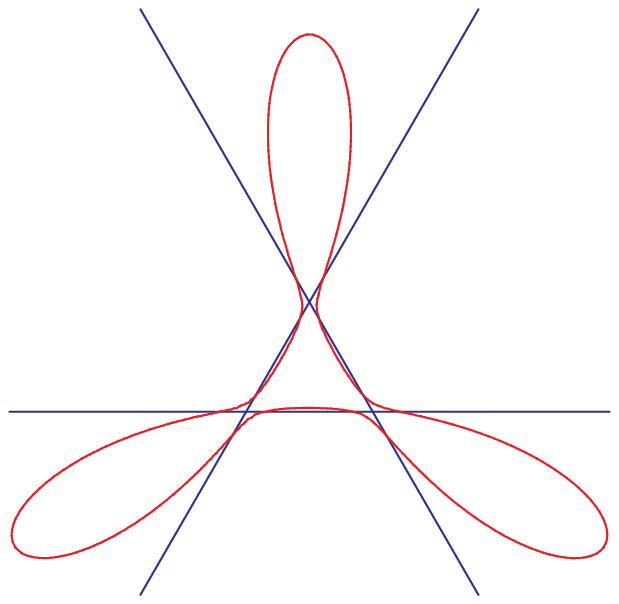}\fi}
%   \put( 52.5,119){$5$}
%   \put( 43, 11){$-2$}
%   \put( 20,51){$-2$}
%   \put( 91,51){$2$} 
   \put(52.5,-6){$(iii)$}
  \end{picture}
   \qquad\qquad
  \begin{picture}(139,145)(0,-8)
   \put(0,-1.8){\ifpictures\includegraphics[height=120pt]{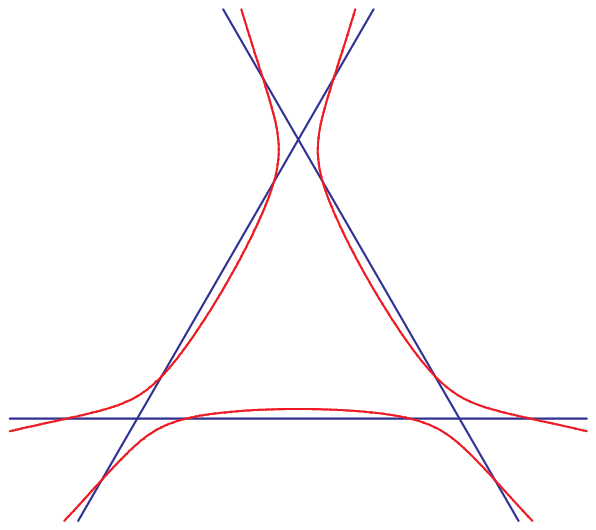}\fi}
%   \put( 55, 87){$1$}
%   \put( 48, 0){$-1$}
%   \put( 13,53){$-1$}
%   \put(107,53){$1$} 
   \put(56.5,-6){$(iii')$}
  \end{picture}
\]
\caption{The cubic and quartic for symmetric configurations.\label{F:cu_qu}}
\end{figure}
%%%%%%%%%%%%%%%%%%%%%%%%%%%%%%%%%%%%%%%%%%%%%%%%%%%%%%%%%%%%%%%%%%%%%%%%%%%%%%%%
 Each point of intersection gives a real tangent, so there are 12 tangents to
 four spheres of equal radii $1.425$ with centers at the vertices of the
 regular tetrahedron with edge length $2\sqrt{2}$.

 One may also see this number 12 using group theory.
 The symmetry group of the tetrahedron, which is the group of permutations of
 the spheres, acts transitively on their common tangents and the isotropy group
 of any tangent has order 2.
 To see this, orient a common tangent and suppose that it meets the spheres
 $a,b,c,d$ in order.
 Then the permutation $(a,d)(b,c)$ fixes that tangent but reverses its orientation,
 and the identity is the only other permutation fixing that tangent.
 \hfill\QED
\end{ex}
%%%%%%%%%%%%%%%%%%%%%%%%%%%%%%%%%%%%%%%%%%%%%%%%%%%%%%%%%%%%%%%%%%%%%%%%%%%%%%%%

This example shows that the bound of 12 common tangents from
Lemma~\ref{L:cubic_quartic} is in fact attained.

\begin{thm}
   There are at most $12$ common real tangent lines to four spheres whose centers are
   not coplanar, and there exist spheres with $12$ common real tangents. 
\end{thm}

\begin{figure}[htb]
\[
  \begin{picture}(200,190)
   \put(0,0){\ifpictures\includegraphics[height=190pt]{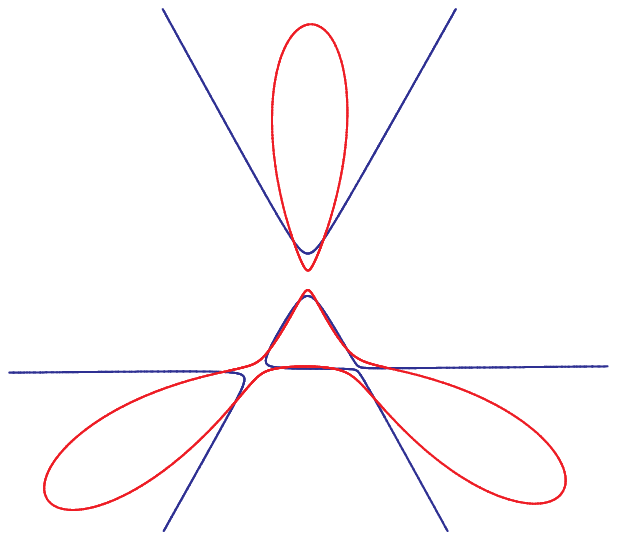}\fi}
%   \put( 80,148){$4$}
%   \put( 80,107){$2$}
%   \put( 80, 27){$-2$}
%   \put(  6,80){$-4$}
%   \put( 48,80){$-2$}
%   \put(138,80){$2$} 
%   \put(177,80){$4$} 
  \end{picture}
   \qquad
  \begin{picture}(200,190)
   \put(0,0){\ifpictures\includegraphics[height=195pt]{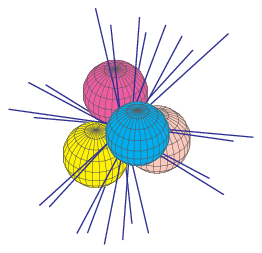}\fi}
  \end{picture}
\]
 \caption{Spheres with 12 common tangents.\label{F:cu_qu_ineq}}
 \end{figure}

%%%%%%%%%%%%%%%%%%%%%%%%%%%%%%%%%%%%%%%%%%%%%%%%%%%%%%%%%%%%%%%%%%%%%%%%%%%%%%%%
\begin{ex}
 We give an example when the radii are distinct,
 namely $1.4$, $1.42$, $1.45$, and $1.474$.
 Figure~\ref{F:cu_qu_ineq} 
 shows the quartic and cubic and the configuration of
 4 spheres and their 12 common tangents.\hfill\QED
\end{ex}
%%%%%%%%%%%%%%%%%%%%%%%%%%%%%%%%%%%%%%%%%%%%%%%%%%%%%%%%%%%%%%%%%
%%%%%%%%%%%%%%%%%%%%%%%%%%%%%%%%%%%%%%%%%%%%%%%%%%%%%%%%%%%%%%%%%%%%%%%%%%%%
%%%%%%%%%%%%%%%%%%%%%%%%%%%%%%%%%%%%%%%%%%%%%%%%%%%%%%%%%%%%%%%%%%%%

\begin{figure}[htb]
\ifpictures
\[
   \includegraphics[height=170pt]{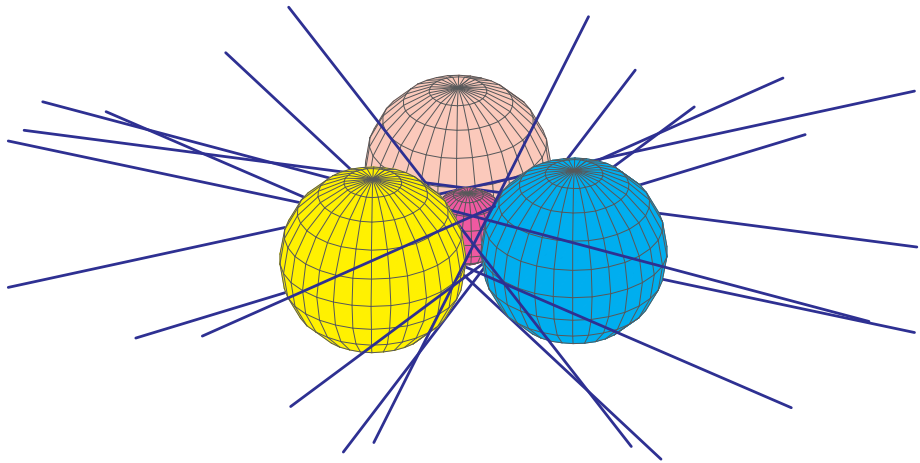}
\]
\fi
\caption{Spheres with coplanar centers and 12 common tangents.\label{F:coplanar}}
\end{figure}

Megyesi investigated configurations of spheres with coplanar centers~\cite{me-01}.
A continuity argument shows that four general such spheres will have 12
complex common tangents (or infinitely many, but this possibility is precluded
by the following example).
Three spheres of radius $4/5$ centered at the vertices of an
equilateral triangle with side length $\sqrt{3}$ and one of radius $1/3$
at the triangle's center have 12 common real tangents.
We display this configuration in Figure~\ref{F:coplanar}.
%%%%%%%%%%%%%%%%%%%%%%%%%%%%%%%%%%%%%%%%%%%%%%%%%%%%%%%%%%%%%%%%%%%%%%%
%%%%%%%%%%%%%%%%%%%%%%%%%%%%%%%%%%%%%%%%%%%%%%%%%%%%%%%%%%%%%%%%%%%%%%%
This configuration of spheres has symmetry group $\Z_2\times D_3$, 
which has order 12 and acts faithfully and transitively on the common tangents. 

Note that these coplanar spheres, unlike the spheres of
Example~\ref{Ex:cubic_quartic}, have unequal radii. 
Megyesi showed that this is no coincidence.

%%%%%%%%%%%%%%%%%%%%%%%%%%%%%%%%%%%%%%%%%%%%%%%%%%%%%%%%%%%%%%%%%%%%%%%%%%%%%%
\begin{thm}\label{Th:Me:8}
 If the centers of four unit spheres in $\R^3$ are coplanar but not colinear,
 then they have at most eight common real tangents, and this bound
 is sharp.
\end{thm}
%%%%%%%%%%%%%%%%%%%%%%%%%%%%%%%%%%%%%%%%%%%%%%%%%%%%%%%%%%%%%%%%%%%%%%%%%%%%%%

The sharpness is shown by the configuration of four spheres of radius $9/10$
with centers at the vertices of a rhombus in the plane $z=0$,
$(0,0,0)^T$, $(2,0,0)^T$, $(1\sqrt{3},0)^T$, and $(3,\sqrt{3},0)^T$,
which have 8 common tangents.
This construction is found in~\cite[Section 4]{theobald-03}.
Observe that the spheres here are disjoint.
%%%%%%%%%%%%%%%%%%%%%%%%%%%%%%%%%%%%%%%%%%%%%%%%%%%%%%%%%%%%%%%%%%%%%%%%%%
\begin{figure}[htb]
\ifpictures
\[
   \includegraphics[height=120pt]{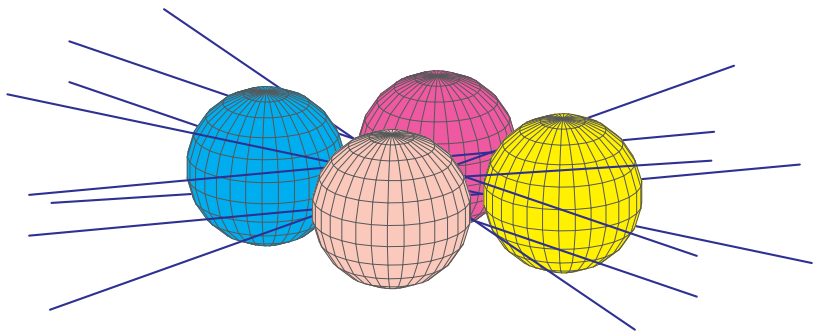}
\]
\fi
\caption{Four unit spheres with coplanar centers and 8 common
     tangents.\label{F:eight}} 
\end{figure}
%%%%%%%%%%%%%%%%%%%%%%%%%%%%%%%%%%%%%%%%%%%%%%%%%%%%%%%%%%%%%%%%%%%%%%%

In the symmetric configuration of Example~\ref{Ex:cubic_quartic} having 12
common tangents, every pair of spheres meet.
%Each pair of spheres meet in the symmetric configuration of
%Example~\ref{Ex:cubic_quartic}.
It is however not necessary for the spheres to meet
pairwise when there are 12 common tangents.
In fact, in both Figures~\ref{F:cu_qu_ineq} and~\ref{F:coplanar}  not all pairs
of spheres meet. 
However, the union of the spheres is connected.
This is not necessary.
In the tetrahedral configuration, if one sphere has radius 1.38 and the other
three have equal radii of 1.44, then the first sphere does not meet the others,
but there are still 12 tangents.

\begin{figure}[htb]
\ifpictures
\[
   \includegraphics[height=150pt]{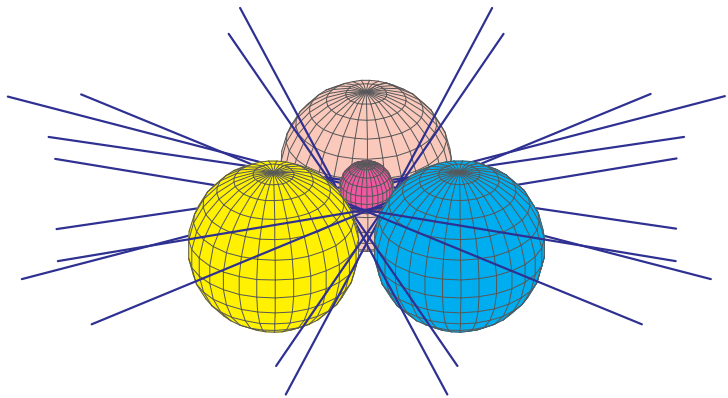}
\]
\fi
\caption{Four disjoint spheres with 12 common tangents.\label{F:disjoint}}
\end{figure}

More interestingly, it is possible to have 12 common real tangents to four
\Red{\emph{disjoint}} spheres.  
Figure~\ref{F:disjoint} displays such a configuration.
%%%%%%%%%%%%%%%%%%%%%%%%%%%%%%%%%%%%%%%%%%%%%%%%%%%%%%%%%%%%%%%%%%%%%%%%
%%%%%%%%%%%%%%%%%%%%%%%%%%%%%%%%%%%%%%%%%%%%%%%%%%%%%%%%%%%%%%%%%%%%%%%%
The three large spheres have radius $4/5$ and are centered at the vertices of an
equilateral triangle of side length $\sqrt{3}$, while the smaller sphere has
radius $1/4$ and is centered on the axis of symmetry of the triangle, but at a
distance of $35/100$ from the plane of the triangle.
It remains an open question whether it is possible for four disjoint
\Red{\emph{unit}} spheres to have 12 common tangents.

%%%%%%%%%%%%%%%%%%%%%%%%%%%%%%%%%%%%%%%%%%%%%%%%%%%%%%%%%%%%%%%%%%%%%%%%%%%%%%
Macdonald, Pach, and Theobald~\cite{mpt-2001} also addressed the question of
degenerate configurations of spheres. 

\begin{thm}
   Four degenerate spheres of equal radii have colinear centers.
\end{thm}

The possible degenerate configurations of spheres having unequal radii remained
open until recently 
when it was settled by Borcea, Goaoc, Lazard, and Petitjean~\cite{bglp-06}.

\begin{thm}
     Four degenerate spheres have colinear centers.
\end{thm}

\begin{figure}[htb]
\[
  \begin{picture}(115,117)(0,-7)
   \put(0, 0){\ifpictures\includegraphics[height=110pt]{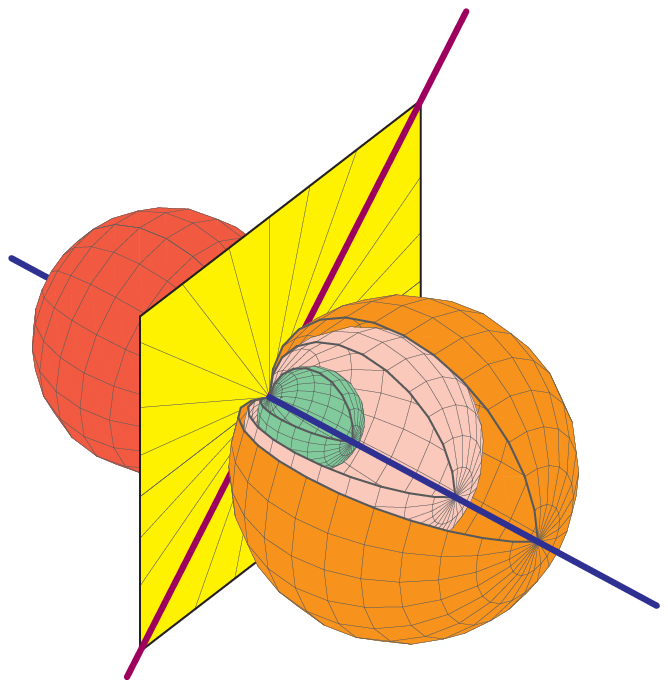}\fi}
   \put(77,100){$\ell$}   \put(90,-1){$m$}
   \put(35,-7){$(i)$}
  \end{picture}
 \qquad
  \begin{picture}(245,115)(0,-5)
   \put(0, 0){\ifpictures\includegraphics[height=110pt]{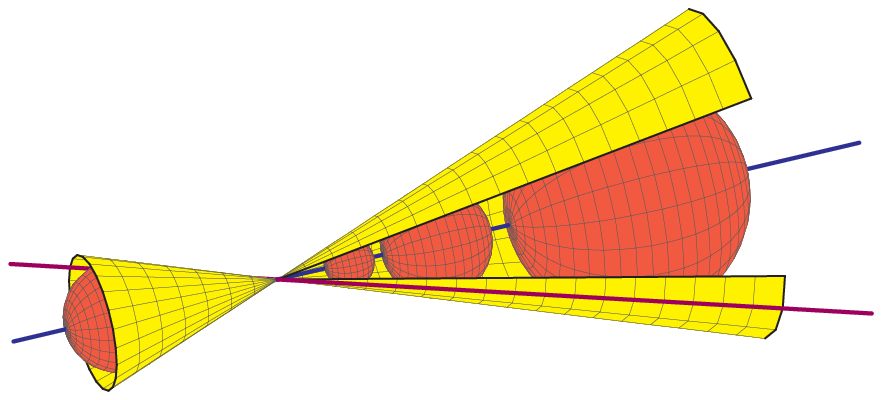}\fi}
   \put(230,10){$\ell$}  \put(220,55){$m$}
   \put(85,-5){$(ii)$}
  \end{picture}
\]

\[
  \begin{picture}(115,115)(0,-12)
   \put(0, -3){\ifpictures\includegraphics[height=110pt]{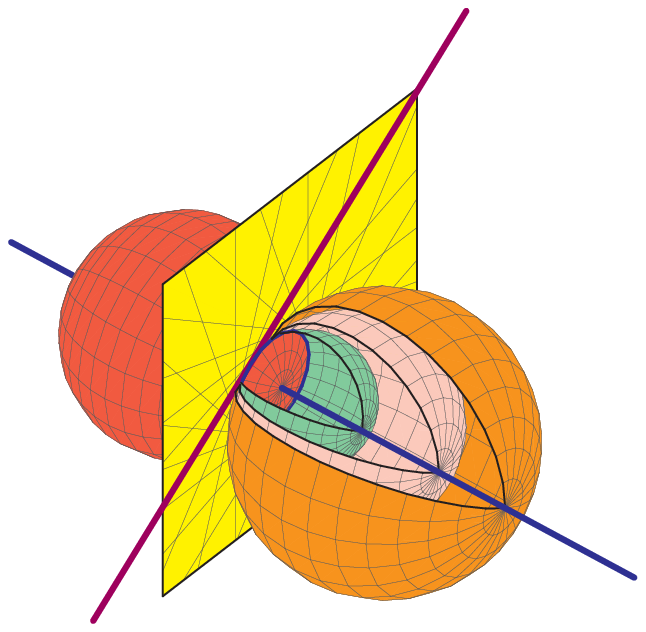}\fi}
   \put(80,94){$\ell$}   \put(95,15){$m$}
   \put(35,-12){$(iii)$}
  \end{picture}
 \qquad
  \begin{picture}(245,115)(0,-12)
   \put(0,5){\ifpictures\includegraphics[width=250pt]{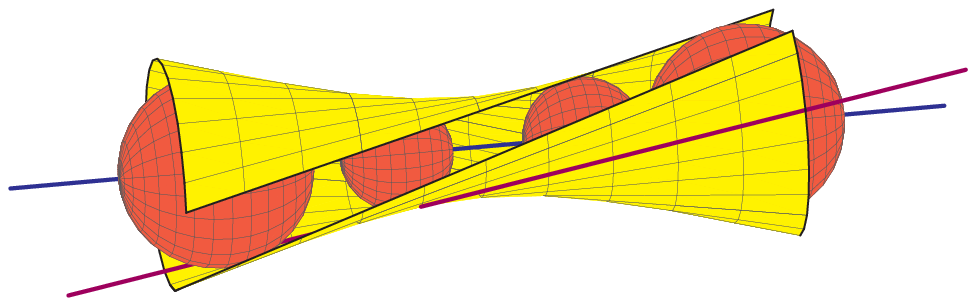}\fi}
   \put(230,35){$m$}
   \put(230,85){$\ell$}
   \put(120,-12){$(iv)$}
  \end{picture}
\]    
\caption{Configurations of four spheres with infinitely many common
   tangents.\label{F:colinear}} 
\end{figure}

Contrary to expectations gained in other work (see Section~\ref{S:degenerate}), their
proof was refreshingly elementary.
Beginning with the formulation of Lemma~\ref{L:cubic_quartic}, they use
algebraic manipulation (with $v^2=0$ playing an important role) to derive a
contradiction to the assumption that the four spheres have infinitely many
common tangents. 
Using a different formulation when the centers of the spheres span a 
plane, they again obtain a contradiction to their being degenerate.
This leaves only the possibility that the four spheres have colinear centers, and 
they then classify degenerate configurations of spheres with colinear centers.
These are displayed in Figure~\ref{F:colinear}.
%%%%%%%%%%%%%%%%%%%%%%%%%%%%%%%%%%%%%%%%%%%%%%%%%%%%%%%%%%%%%%%%%%%%%%%%%
%%%%%%%%%%%%%%%%%%%%%%%%%%%%%%%%%%%%%%%%%%%%%%%%%%%%%%%%%%%%%%%%%%%%%%%%%%%%%%%%
The spheres  $(i)$ are tangent to each other at a point, or $(ii)$ 
are inscribed in a circular cone or cylinder, or $(iii)$ meet in a common circle,
or $(iv)$ are inscribed in a hyperboloid 
of revolution. 
The common tangents are generated from a single common tangent line $\ell$ by
rotation about the line $m$ containing their centers.

%%%%%%%%%%%%%%%%%%%%%%%%%%%%%%%%%%%%%%%%%%%%%%%%%%%%%%%%%%%%%%%%%%%%%
\subsection{Common tangents to $k$ spheres which meet $4{-}k$ lines}
 Theobald~\cite{theobald-enumerative} studied the intermediate cases of our
 algebraic core problem, those involving {\it both} lines and spheres.

%%%%%%%%%%%%%%%%%%%%%%%%%%%%%%%%%%%
\begin{thm}
 Bounds for the number $N_k$ of common tangents to $k$ spheres which also meet
 $4-k$ fixed lines for $1\leq k\leq 3$ are as follows
\[
   N_1\ \leq\ 4\,,\qquad
   N_2\ \leq\ 8\,,\quad\mbox{and}\quad
   N_3\ \leq\ \Blue{12}\,,
\]
 and these are sharp.
 There exist unit spheres and lines having exactly these numbers of real common
 tangents. 
\end{thm}
%%%%%%%%%%%%%%%%%%%%%%%%%%%%%%%%%%%

\noindent{\it Proof of bounds.}
 In Pl\"ucker coordinates, the set of transversals and tangents is defined by 
 $4{-}k$ linear equations and $k{+}1$ quadratic equations.
 Indeed, the condition for a line to meet a fixed line~\eqref{Eq:linear_Pluecker}
 is linear in the Pl\"ucker coordinates, the condition for a line to be
 tangent to a sphere~\eqref{Eq:Sphere_eq} is quadratic in the Pl\"ucker
 coordinates, and the set of lines is defined by the quadratic Pl\"ucker
 equation~\eqref{Eq:Pluecker}.
 B\'ezout's Theorem gives the bounds
\[
   N_1\ \leq\ 4\,,\qquad
   N_2\ \leq\ 8\,,\quad\mbox{and}\quad
   N_3\ \leq\ \Blue{16}\,.
\]

 The B\'ezout bound of 16 for $N_3$ includes not only the geometrically
 meaningful tangents to the three spheres which meet a line, $\ell$, but also a
 contribution from the lines tangent to the spherical
 conic~\eqref{Eq:Spherical_Conic} at infinity which also meet $\ell$.
 Since there are two such lines, the contribution is $2m$, where $m$ is the
 algebraic multiplicity of either line in the intersection.
 Theobald~\cite{theobald-enumerative} computes $m=2$, and so there are at most
 12 lines tangent to three spheres which also meet a line.
\QED\medskip

We now give constructions which realize these bounds.

%%%%%%%%%%%%%%%%%%%%%%%%%%%%%%%%%%%%%%%%%%%%%%%%%%%%%%%%%%%%%%%%%%%%%%%%%%%%%%
\subsubsection{One sphere and three lines}
Suppose that the three lines are in the (slightly) degenerate position where
$\ell_1$ and $\ell_2$ are skew, but each meets $\ell_3$ in points 
$p_1:=\ell_2\cap\ell_3$ and $p_2:=\ell_1\cap\ell_3$.
Then the common transversals to the three lines form two pencils,
the lines through $p_1$ which meet $\ell_1$ and the lines through $p_2$ which
meet $\ell_2$.
Let $H_1$ be the plane spanned by the first pencil and $H_2$ the plane spanned
by the second pencil. 

\begin{figure}[htb]
\[
  \begin{picture}(255,165)
   \put(0,0){\includegraphics[height=165pt]{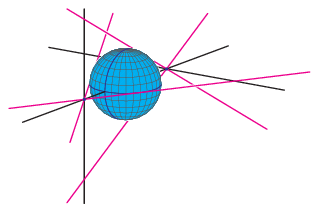}}
%   \put(0,0){\includegraphics[height=165pt]{figures/shell.eps}}
   \put( 71, 1){$\ell_1$}
   \put(228,85){$\ell_2$}
   \put( 12,57){$\ell_3$}

   \put( 51, 95){$p_2$}  \put(132,128){$p_1$}
   \put( 94,157){$C_2$}  \put( 95,154){\vector(-1,-3){18}}
   \put(141, 80){$C_1$}  \put(138,85){\vector(-1,0){50}}

   \put(-10, 87){$t_{21}$}   \put( 47,42){$t_{22}$}
   \put(156,162){$t_{11}$}   \put(214,52){$t_{12}$}

  \end{picture}
\]
\caption{One sphere and three lines with four common tangents.\label{F:3l1S}}
\end{figure}

 If $S$ is a sphere which meets both $H_1$ and $H_2$ in circles 
$C_1$ and $C_2$ and contains neither $p_1$ nor $p_2$, then four of the
transversals to the three 
lines will be tangent to $S$. 
Indeed, in each plane $H_i$ there are two tangents $t_{i1}$ and $t_{i2}$ to the
circle  $C_i$ through $p_i$, as $p_i$ is exterior to $C_i$.
See Figure~\ref{F:3l1S} for a picture of such a configuration.
There, the line $\ell_3$ is the $x$-axis, $\ell_1$ contains the point
$p_2:=(2,0,0)$ and is parallel to the $z$-axis, $\ell_2$ contains the point
$p_1:=(-2,0,0)$ and is parallel to the $y$-axis, and the sphere is centered at the
origin with radius $1$. \hfill\QED
%%%%%%%%%%%%%%%%%%%%%%%%%%%%%%%%%%%%%%%%%%%%%%%%%%%%%%%%%%%%%%%%%%%%%%%%%%
%%%%%%%%%%%%%%%%%%%%%%%%%%%%%%%%%%%%%%%%%%%%%%%%%%%%%%%%%%%%%%%%%%%%%%%%%%

%%%%%%%%%%%%%%%%%%%%%%%%%%%%%%%%%%%%%%%%%%%%%%%%%%%%%%%%%%%%%%%%%%%%%%%%%%%%%%
\subsubsection{Two spheres and two lines}
Suppose the lines $\ell_1$ and $\ell_2$ meet in a point $p$ and span a plane $H$.
Then there are two families of transversals to the two lines;
lines that pass through $p$ and lines contained in $H$.
If $H$ meets the spheres in disjoint circles, then the four common
tangents to the circles in $H$ are common tangents to the spheres which meet 
both $\ell_1$ and $\ell_2$.
Dually, the point $p$ and spheres can be located so that there are four common
tangents to the spheres through $p$.

Supppose that $S_1$ and $S_2$ are spheres having radius $\sqrt{5}$
and centers $(\pm2,0,0)$.
If $\ell_1$ and $\ell_2$ have parametrization
$(t,2\mp2t,\pm3t)$ for $t\in\R$, then they meet in the point $p=(0,2,0)$ and
span the plane $H$ defined by $3y+2z=6$.
We claim that there are eight common transversals to $\ell_1$ and $\ell_2$
which are tangent to $S_1$ and $S_2$.

\begin{figure}[htb]
 \[
   \begin{picture}(175,152)
   \put(0,0){\includegraphics[height=150pt]{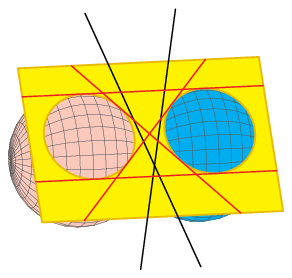}}

   \put( 37,142){$\ell_1$}   \put(103,142){$\ell_2$}
   \put(  0,120){$H$}   

   \put( 33,124){$t_1$}   \put(117,128){$t_2$}
   \put(155,103){$t_3$}   \put(163, 55){$t_4$}

  \end{picture} 
\qquad\qquad
  \begin{picture}(190,152)
   \put(0,0){\includegraphics[height=150pt]{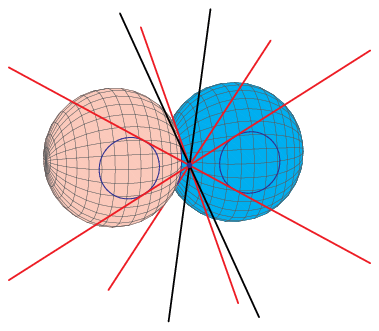}}
   \put( 48,142){$\ell_1$}   \put(107,142){$\ell_2$}
   \put( 4, 27){$t_5$}      \put( 11,123){$t_6$}
   \put(74,136){$t_7$}      \put(136,135){$t_8$}
  \end{picture} 
\]
\caption{Eight tangents to two spheres that meet two lines.\label{F:four_four}}
\end{figure}

Indeed, the plane $H$ meets the spheres in disjoint circles, contributing four
common tangents $t_1,\dotsc,t_4$.
There are four common tangents $t_5,\dotsc,t_8$ through the point $p$;
two in the plane $y=2$ and two lying in the $yz$-plane.
These have parametrizations
 \[
   t_5,t_6\ =\ (t, 2, \pm t/\sqrt{3})\qquad\mbox{and}\qquad
   t_7,t_8\ =\ (0, 2+t, \pm t/\sqrt{3})\,.
 \]
We display this configuration in Figure~\ref{F:four_four}, showing the common
tangents in the plane $H$ on the left and those through the point $p$ on the right.
In the second picture, we draw the circles where the plane $y=2$ meets the
spheres. \hfill\QED
%%%%%%%%%%%%%%%%%%%%%%%%%%%%%%%%%%%%%%%%%%%%%%%%%%%
%%%%%%%%%%%%%%%%%%%%%%%%%%%%%%%%%%%%%%%%%%%%%%%%%%%%%%%%%%%%%%%%%%%%%%%%%%

%%%%%%%%%%%%%%%%%%%%%%%%%%%%%%%%%%%%%%%%%%%%%%%%%%%%%%%%%%%%%%%%%%%%%%%%%%%%%%
\subsubsection{Three spheres and one line}
Suppose that $S_1$, $S_2$, and $S_3$ are unit spheres centered at the vertices of
an equilateral triangle with edge length $e$, where $\sqrt{3}<e<2$, and
that $\ell$ is the line perpendicular to the plane of the triangle passing
through its center.
Then there are 12 common real tangents to the sphere which meet $\ell$.

\begin{figure}[htb]
 \[
   \begin{picture}(190,147)(0,-3)
   \put(0,0){\includegraphics[height=150pt]{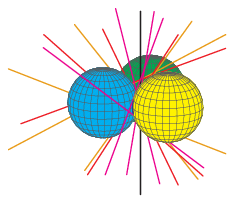}}
   \put(99,-5){$\ell$}
   \end{picture}
 \]
 \caption{Three spheres and one line with 12 common tangents.\label{F:1l3s}}
\end{figure}

To see this, note that any common tangent lies in a vertical plane
through $\ell$ which meets all three spheres. 
Such a plane meets each sphere in a circle, and one of the three circles
necessarily contains another, unless the plane contains a vertex of the triangle.
A plane through $\ell$ and through a vertex meets the spheres in two disjoint
circles---one is the intersection of two spheres.
The four common tangents to these circles in this plane are common tangents
to the spheres which meet $\ell$.
As there are three such planes, we obtain 12 common tangents.
Figure~\ref{F:1l3s} shows this configuration when $e=1.9$.\hfill\QED
%%%%%%%%%%%%%%%%%%%%%%%%%%%%%%%%%%%%%%%%%%%%%%%%%%%%%%%%%%%%%%%%%%%%%%
%%%%%%%%%%%%%%%%%%%%%%%%%%%%%%%%%%%%%%%%%%%%%%%%%%%%%%%%%%%%%%%%%%%%%

%%%%%%%%%%%%%%%%%%%%%%%%%%%%%%%%%%%%%%%%%%%%%%%%%%%%%%%%%%%%%%%%%%%%%
%
%
\subsection{Degenerate configurations}\label{S:degenerate}

We characterize configurations of $k$ spheres and $4{-}k$ lines which are
degenerate in that they have 
%infinitely many
a 1-dimensional family of 
%%%%%%%%%%%%%%%%%%%%%%%%%%%%%
common tangents.

%%%%%%%%%%%%%%%%%%%%%%%%%%%%%%%%%%%%%%%%%%%%%%%%%%%%%%%%%%%%%%%%%%%%%
%
%
\subsubsection{One sphere and three lines}
We leave the case when some of the lines meet as an exercise to the reader, and
suppose that the three lines are pairwise skew so that their common transversals
are one of the families of lines on a hyperboloid $H$ which contains them as in
Figure~\ref{F:3lines}. 
Each point where the sphere is tangent to $H$ gives a common transversal which
is tangent to the sphere.
Thus the configuration of lines and the sphere is degenerate if and only if the
sphere is tangent to $H$ along a curve, $C$.

  The intersection of the sphere with the hyperboloid $H$ will have multiplicity
at least~2 along $C$.
Since $C$ cannot be a line, and B\'ezout's Theorem implies that the product of
this multiplicity and the degree of $C$ is at most 4, we conclude that the
multiplicity is 2 and that $C$ has degree 2.
Thus $C$ is a plane conic lying on the sphere, and is thus a circle.
This will imply that the hyperboloid is a right circular hyperboloid inscribing
the sphere.
We do not give a picture as it is similar to Figures~\ref{F:colinear}({\it iv})
and~\ref{F:Thm1}({\it iv}).  

%%%%%%%%%%%%%%%%%%%%%%%%%%%%%%%%%%%%%%%%%%%%%%%%%%%%%%%%%%%%%%%%%%%%%
%
%
\subsubsection{Two spheres and two lines}
The most interesting geometry occurs when we study two spheres and two lines.
Megyesi, Sottile, and Theobald~\cite{mst-2003} considered the more general
situation where we replace the spheres by quadrics. 
Fix two skew lines $\ell_1$ and $\ell_2$ and a quadric $Q$ in general position.
Then the set $\tau$ of tangents to $Q$ which meet both $\ell_1$ and $\ell_2$ is
an irreducible algebraic curve (in fact it has genus 1).
We consider the set $\calQ$ of quadrics $Q$ which are tangent to every line in
$\tau$.
This set $\calQ$ consists of quadrics $Q'$ such that $Q'$, $Q$, 
$\ell_1$, and $\ell_2$ are degenerate.

Since a quadric is defined as the zero set of a quadric form, quadrics
correspond to symmetric $4\times 4$ matrices, modulo scalars.
This identifies $\P^9$ as the space of quadrics.
The following theorem describes the set $\calQ$ as a complex algebraic
variety. 

\begin{thm}\label{thm:24-12}
 The set $\calQ$ of quadrics having degenerate position with respect to skew
 lines $\ell_1$  and $\ell_2$ and a general quadric~$Q$ in $\mathbb{P}^3$
 is an algebraic curve in $\P^9$ of degree $24$  having $12$ irreducible
 components, each a plane conic.
\end{thm}

The proof of this statement uses some geometric reductions, and then a very
challenging symbolic computation using the computer algebra system
{\sc Singular}~\cite{GPS05}, which uses Gr\"obner bases to manipulate
polynomials.

That $\calQ$ is composed of algebraic curves is not too hard to see and it
involves some beautiful geometry. 
It turns out that the set of non-zero real numbers $\R^\times$ (or $\C^\times$
if we are working over the complex numbers) acts on $\P^3$ with fixed points
$\ell_1$ 
and $\ell_2$ so that the orbits are the transversals to $\ell_1$ and $\ell_2$.
One way to see this is to choose coordinates so that $\ell_1$ is the $x$-axis
and $\ell_2$ is the line at infinity in the $yz$-plane.
Then the common transversals to $\ell_1$ and $\ell_2$ are the lines
perpendicular to the $x$-axis, and $\R^\times$ acts by simultaneously scaling
the $y$- and $z$-coordinates.
It follows that if $Q'$ is tangent to every line of $\tau$, then so is any
translate of $Q'$. 

Not all of the 12 families of quadrics in $\calQ$ will contain real smooth
quadrics.
From the calculations in~\cite{mst-2003}, which are archived on a
webpage\footnote{{\tt http://www.math.tamu.edu/\~{}sottile/pages/2l2s/index.html}}, 
if $Q$ is a sphere or ellipsoid then four of the families will be real, and
there is a choice of $Q$ for which all 12 families are real.
Furter analysis reveals that if $\ell_1$ and $\ell_2$ are in $\R^3$ and if $Q$
is a sphere, then there are no other spheres in $\calQ$.

\begin{figure}[htb]
\[
  \begin{picture}(375,115)(-10,5)
   \put(10,15){\ifpictures\includegraphics[height=95pt]{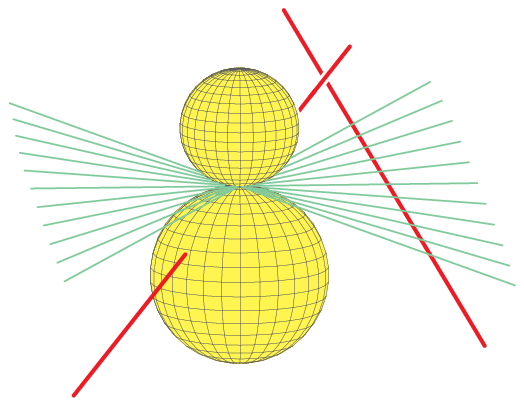}\fi}
    \put(-10,55){$\tau$}
    \put(  1, 4){$\ell_1$}  \put(119,15){$\ell_2$}
    \put(52,5){({\it i})} 

   \put(160,25){\ifpictures\includegraphics[height=70pt]{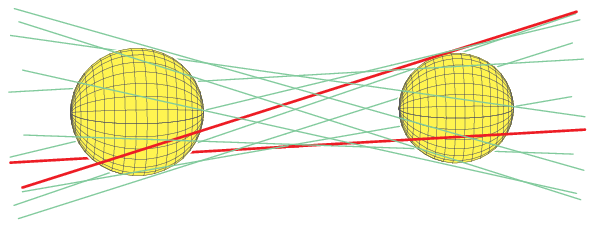}\fi}
   \put(330,25){$\tau$}
   \put(355,93){$\ell_1$}  \put(355,53){$\ell_2$}
   \put(255,5){({\it ii})}
  \end{picture}
\]
\caption{Examples from Theorem~\ref{th:2l2sclassification}.\label{fig3}}
\end{figure}

Theorem~\ref{thm:24-12} underlies the following classification of degenerate
configurations from~\cite{mst-2003}, which is illustrated in
Figure~\ref{fig3}. 

\begin{thm}
\label{th:2l2sclassification} 
Let $S_1 \neq S_2$ be spheres in $\R^3$,
and let $\ell_1 \neq \ell_2$ be two skew lines in $\R^3$.
There are infinitely many transversals to $\ell_1$ and $\ell_2$ which are tangent to
$S_1$ and $S_2$ in exactly the following cases.
\begin{enumerate}
 \item[({\it i})] $S_1$ and $S_2$ are tangent to each other at
            a point $p$ which lies on one line, and the second line lies in
            the common tangent plane to the spheres at the point $p$.
 \item[({\it ii})] The lines $\ell_1$ and $\ell_2$ are each tangent to both $S_1$
            and $S_2$ and they are images of each other under a
            rotation about the line connecting the centers of $S_1$ and $S_2$.
\end{enumerate}
\end{thm}

%%%%%%%%%%%%%%%%%%%%%%%%%%%%%%%%%%%%%%%%%%%%%%%%%%%%%%%%%%%%
%%%%%%%%%%%%%%%%%%%%%%%%%%%%%%%%%%%%%%%%%%%%%%%%%%%%%%%%%%%%

%%%%%%%%%%%%%%%%%%%%%%%%%%%%%%%%%%%%%%%%%%%%%%%%%%%%%%%%%%%%%%%%%%%%%
%
%
\subsubsection{Three spheres and one line}
\begin{figure}[htb]
\[
  \begin{picture}(128,110)(-20,-5)
   \put(0,10){\ifpictures\includegraphics[height=100pt]{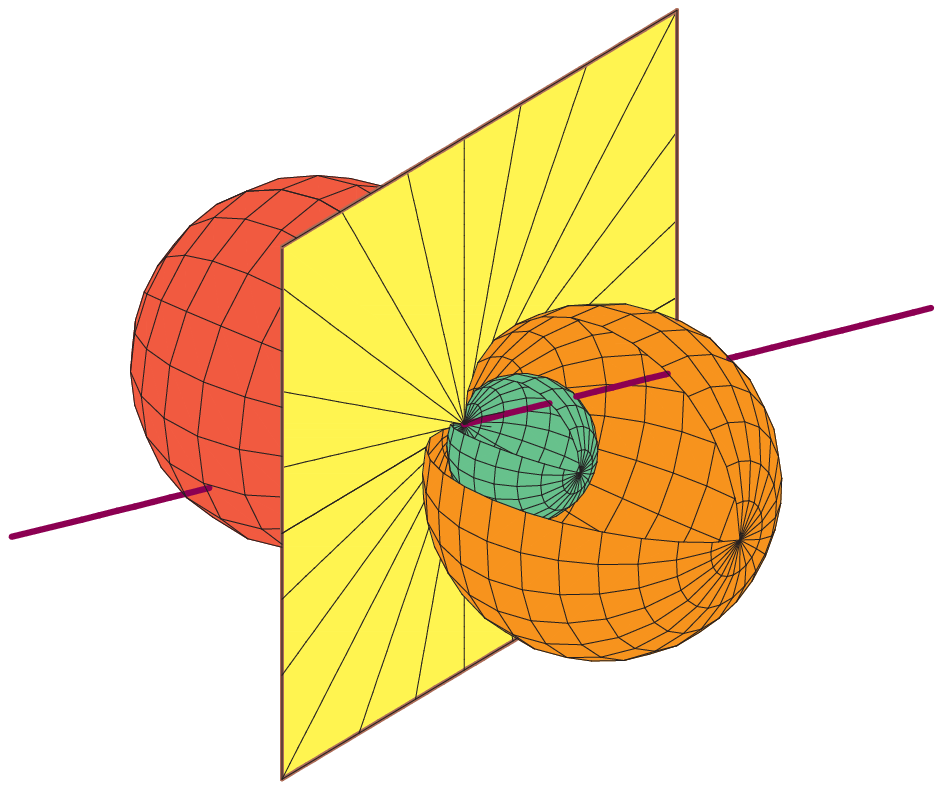}\fi}
   \put(100,74){$\ell$}
   \put(40,-5){({\it ia})}
  \end{picture}
 \qquad
  \begin{picture}(85,110)(0,-5)
   \put(0,10){\ifpictures\includegraphics[height=100pt]{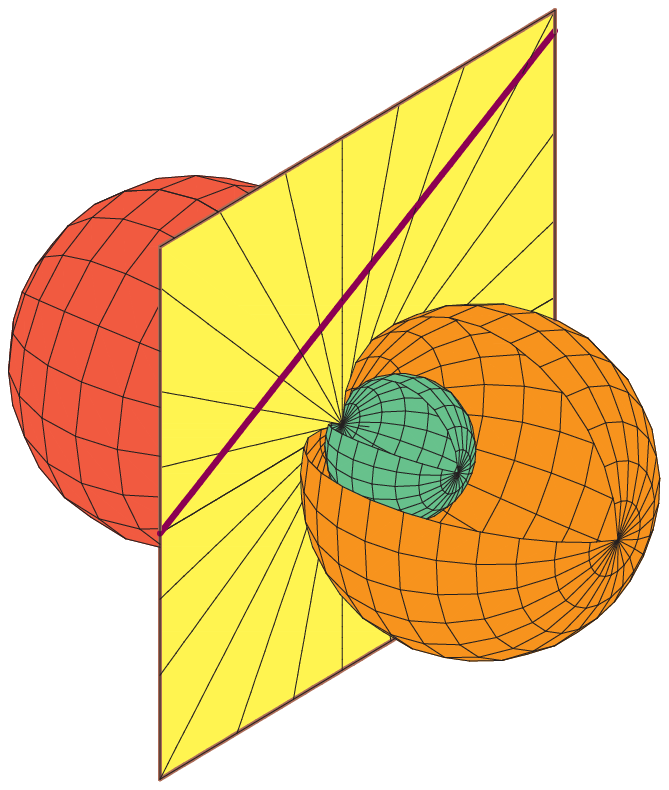}\fi}
   \put(75,105){$\ell$}
   \put(35,-5){({\it ib})}
   \end{picture}
 \qquad
  \begin{picture}(175,110)(0,-5)
   \put(0,10){\ifpictures\includegraphics[height=100pt]{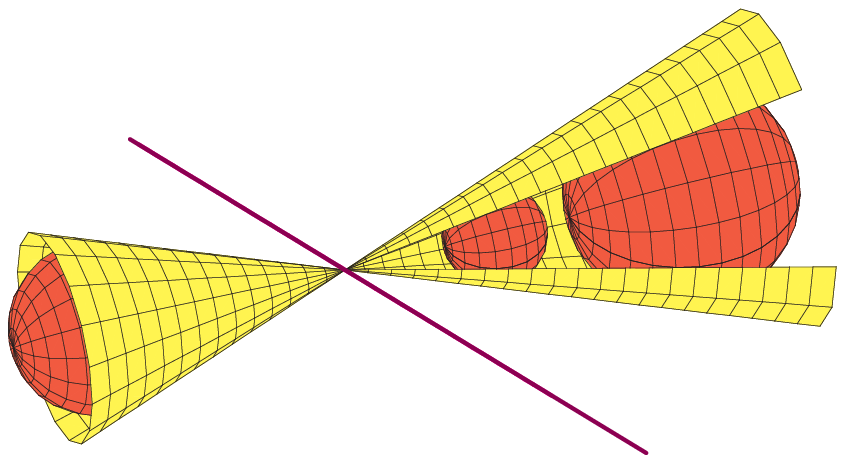}\fi}
   \put(40,77){$\ell$}
   \put(85,-5){({\it ii})}
  \end{picture}
\]

\[
  \begin{picture}(105,110)(0,-10)
   \put(0,10){\ifpictures\includegraphics[height=100pt]{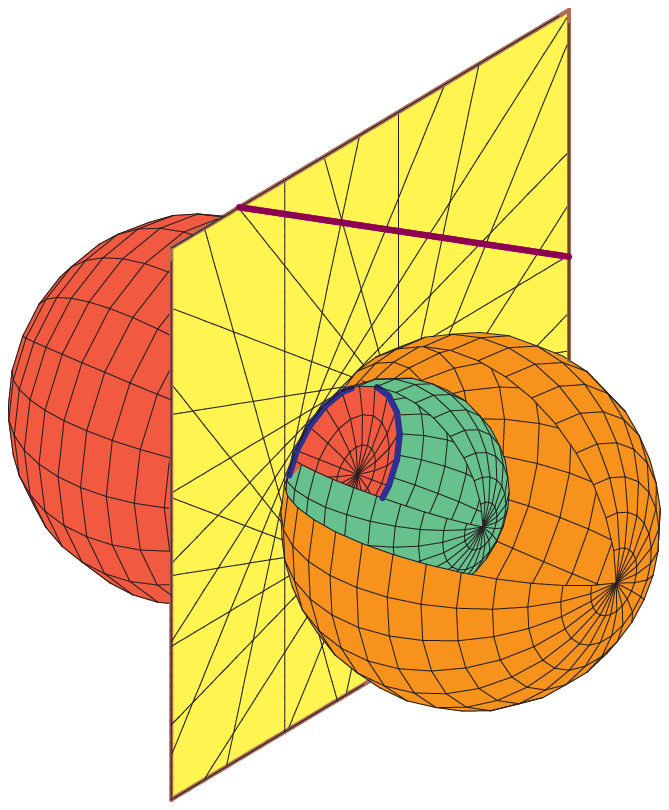}\fi}
   \put(75,80){$\ell$}
   \put(50,-5){({\it iii})}
  \end{picture}
 \qquad
  \begin{picture}(249,83)(0,-10)
   \put(0,15){\ifpictures\includegraphics[height=80pt]{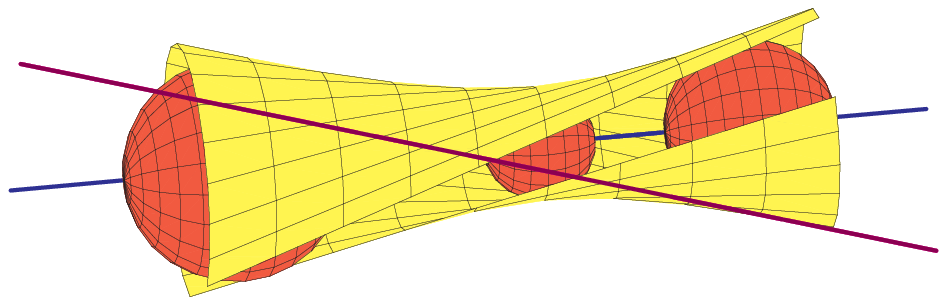}\fi}
   \put(5,35){$m$}
   \put(15,82){$\ell$}
   \put(120,-5){({\it iv})}
  \end{picture}
\]

\caption{The geometry of Theorem~\ref{T:ls2}.\label{F:Thm1}}

\end{figure}

Degenerate positions of three spheres and one line were classified
by Megyesi and Sottile~\cite{megyesi-sottile-2005}.
This was more complicated than the classifications with fewer than three spheres
and did not contain as much interesting geometry as we just saw.
Here is the classification, which is illustrated in Figure~\ref{F:Thm1}.

\begin{thm}\label{T:ls2}
 Let $\ell$ be a line and $S_1$, $S_2$, and $S_3$ be spheres in $\R^3$. 
 Then there are infinitely many lines that meet $\ell$ and are 
 tangent to each sphere in precisely the following cases.
\begin{enumerate}
 \item [({\it i})] The spheres are tangent to each other at the same point 
             and either $($a$)$  $\ell$ meets that point,  
             or $($b$)$ it lies in the common tangent plane, or both. 
             The common tangent lines are the lines in the tangent plane
             meeting the point of tangency.
 \item [({\it ii})] The spheres are tangent to a cone whose apex lies on 
             $\ell$.
             The common tangent lines are the ruling of the cone.

 \item [({\it iii})] The spheres meet in a common circle and the line
             $\ell$ lies in the plane of that circle. 
             The common tangents are the lines in that plane tangent to the
             circle.
 \item [({\it iv})] The centers of the spheres lie on a line  
             $m$ and $\ell$ is tangent to all three spheres.
             The common tangent lines are one ruling on the hyperboloid of
             revolution obtained by rotating $\ell$ about $m$. 
\end{enumerate}
 In cases  $($ib$)$, $($iii$)$ and $($iv$)$ lines parallel to $\ell$ are
 excluded. 
\end{thm}

%%%%%%%%%%%%%%%%%%%%%%%%%%%%%%%%%%%%%%%%%%%%%%%%%%%%%%%%%%%%%%%%%%%%%%%%%
%%%%%%%%%%%%%%%%%%%%%%%%%%%%%%%%%%%%%%%%%%%%%%%%%%%%%%%%%%%%%%%%%%%%%%%%%%%%5

For this, they begin with two spheres and one line, and consider the common
tangents to the spheres which also meet the line. 
These tangents form a curve $\tau$ of degree 8 in Pl\"ucker space.  
When the configuration of the two spheres and the fixed line becomes degenerate, 
this curve becomes reducible.
It is impossible for any component $\sigma$ of $\tau$ to have
degree 3 or 5 or 7, and if a component $\sigma$ has degree 4 or 6 or 8, then it
determines the two spheres.
Only in the cases when $\sigma$ has degree 1 or 2 can there be more than
2 spheres; these are described in Theorem~\ref{T:ls2}.

%%%%%%%%%%%%%%%%%%%%%%%%%%%%%%%%%%%%%%%%%%%%%%%%%%%%%%%%%%%%%%%%%%%%%%%%%%%%5
%
%
%
\section{Open problems}\label{S:openproblems}

As we have seen, the tangent problems lead to
a variety of problems combining discrete and algebraic geometry.
Some open questions have already been mentioned, in particular determining
the maximum number of real common tangent lines to four disjoint unit balls in 
$\R^3$. A related open problem is to determine 
the maximum number of different geometric permutations there can be 
to four unit balls in 3-space. 
It is known that 2 is a lower bound (Figure~\ref{F:eight} is one such
configuration) and 3 is an upper bound (see~\cite{cgh-2005}).
If the balls have different radii, then Figure~\ref{F:disjoint} gives four balls
with six geometric permutations.

We already noted that if we generalize the problems from
spheres to quadric surfaces, the number of solutions can increase.
In particular, four general quadric surfaces in $\R^3$ have 32 complex
common tangent lines, and all these tangent lines can be
real~\cite{sottile-theobald-2005}. 
However, a complete  characterization of the degenerate configurations is not
known~\cite{BGLP2}. 

From the viewpoint of the inherent algebraic degree, the problem becomes more  
difficult for tangents to spheres in $\R^d$.
If $d\geq 3$, then there are $3 \cdot 2^{d-1}$
complex common tangent lines to $2d-2$ general spheres in $\R^d$ and 
all of the tangent lines can be real~\cite{sottile-theobald-2002}.
A characterization of
the degenerate configurations seems currently to be out of reach.
%%%%%%%%%%%%%%%%%%%%%%%%%%%%%%%%%%%%%%%%%%%%%%%%%%%%%%%%%%%%%%%%%%%
%%%%%%%%%%%%%%%%%%%%%%%%%%%%%%%%%%%%%%%%%%%%%%%%%%%%%%%%%%%%%%%%%%%
\providecommand{\bysame}{\leavevmode\hbox to3em{\hrulefill}\thinspace}
\providecommand{\MR}{\relax\ifhmode\unskip\space\fi MR }
% \MRhref is called by the amsart/book/proc definition of \MR.
\providecommand{\MRhref}[2]{%
  \href{http://www.ams.org/mathscinet-getitem?mr=#1}{#2}
}
\providecommand{\href}[2]{#2}

\end{document}
%%%%%%%%%%%%%%%%%%%%%%%%%%%%%%%%%%%%%%%%%%%%%%%%%%%%%%%%%%%%%%%%%%%
%%%%%%%%%%%%%%%%%%%%%%%%%%%%%%%%%%%%%%%%%%%%%%%%%%%%%%%%%%%%%%%%%%%
\bibliographystyle{amsplain}
\bibliography{Snowbird}
\end{document}